\documentclass[10pt]{article}
\usepackage{amsfonts}
\usepackage{stmaryrd}
\usepackage{mathrsfs}
\usepackage{amssymb,amsfonts,amsmath, psfrag,eepic,colordvi,epsfig}
\usepackage{enumerate}
\usepackage{color}
\usepackage{graphicx}
\usepackage{subfigure}
\usepackage{array}
\usepackage{caption}
\usepackage{cite}
\usepackage{hyperref}
\topskip 
\numberwithin{equation}{section}
\newtheorem{theorem}{Theorem}[section]

\newtheorem{corollary}[theorem]{Corollary}

\newtheorem{lemma}[theorem]{Lemma}

\numberwithin{figure}{section}

\begin{document}
\parskip 7pt

\pagenumbering{arabic}
\def\sof{\hfill\rule{2mm}{2mm}}
\def\ls{\leq}
\def\gs{\geq}
\def\SS{\mathcal S}
\def\qq{{\bold q}}
\def\MM{\mathcal M}
\def\TT{\mathcal T}
\def\EE{\mathcal E}
\def\lsp{\mbox{lsp}}
\def\rsp{\mbox{rsp}}
\def\pf{\noindent {\it Proof.} }
\def\mp{\mbox{pyramid}}
\def\mb{\mbox{block}}
\def\mc{\mbox{cross}}
\def\qed{\hfill \rule{4pt}{7pt}}
\def\pf{\noindent {\it Proof.} }
\textheight=22cm

\begin{center}
{\Large\bf Congruences modulo
 arbitrary powers of $5$ and $7$ for
   Andrews and Paule's partition
   diamonds with $(n+1)$ copies of $n$
 }
\end{center}

\begin{center}

Julia Q.D. Du$^{1}$ and Olivia X.M. Yao$^{2}$

$^{1}$School of Mathematical Sciences,\\
Hebei Key Laboratory of Computational Mathematics and Applications,\\
Hebei Research Center of the Basic Discipline Pure Mathematics,\\
Hebei Normal University,
 Shijiazhuang 050024, P.R. China

$^{2}$School of Mathematical Sciences, \\
  Suzhou University of Science and
Technology, \\
 Suzhou,  215009, Jiangsu Province,
 P.R. China

e-mail: qddu@hebtu.edu.cn,
 yaoxiangmei@163.com

 \end{center}

\noindent {\bf Abstract.} Recently,
  Andrews and Paule introduced
   a partition function
    $PDN1(N)$ which
  denotes  the number of partition
   diamonds with $(n+1)$ copies of $n$
    where summing the parts
     at the links gives $N$.
 They also presented
        the generating function for  $PDN1(n)$
         and proved
          several
           congruences modulo
          5,7,25,49 for $PDN1(n)$.
           At the end of their
            paper, Andrews and Paule
             asked for determining
               infinite families of
                congruences similar to
                Ramanujan's classical $
                 p(5^kn +d_k) \equiv
                  0 \pmod {5^k}$, where
             $24d_k\equiv 1 \pmod {5^k}$ and $k\geq 1$.
 In this paper, we give an answer
  of Andrews and Paule's open problem
   by proving  three congruences modulo
 arbitrary powers of
 $5$  for $PDN1(n)$.
 In addition, we prove two congruences modulo
 arbitrary powers of
 $7$  for $PDN1(n)$, which are analogous to Watson's
   congruences   for $p(n)$.



\section{Introduction}

\allowdisplaybreaks

This paper is devoted to
 the study of congruences
  modulo
 arbitrary powers of $5$ and $7$
for the number of
   Andrews and Paule's  partition
   diamonds with $(n+1)$ copies of $n$
    where summing the parts at the links gives $N$.
    Let us begin with
  a class of combinatorial objects: partitions with $n
   $ copies of $n$, which play an auxiliary role in the theory of plane
partitions. Let
  $M$ be the set   of
 subscripted positive integers wherein the subscript
  does not exceed
the integer. We
   order this set lexicographically:
\[
1_1<2_1<2_2<3_1<3_2<3_3<4_1<4_2<\cdots.
\]

We say $\pi=(\pi_1,\pi_2,\ldots, \pi_k)$ is a partition of $m$ with
$n$ copies of $n
 $ if $\pi_1+\pi_2+\cdots+\pi_k=m$
  and $\pi_1\geq \pi_2\geq \cdots
 \geq \pi_k$ with $\pi_i\in M$.
  The study of partitions with $n$
   copies of $n$
    had its origins in Regime III of
    the hard hexagon model \cite{Andrews-Baxter}.
      It was made explicit in \cite{Andrews}
       where two interesting
        results
          were posed
           linking partitions with $n$
            copies of $n$ to
  ordinary partitions.
  For more details on partitions
 with $n$ copies of $n$, one can
  see \cite{Agarwal,Andrews}.

 In addition,  for
  partitions with $ (n + 1) $ copies of $n$,   the
subscript 0 will be allowed, and lexicographic order is maintained,
 namely,
\[
1_0<1_1< 2_0<2_1<2_2<3_0<3_1<3_2<3_3<4_0 <4_1<4_2<\cdots.
\]
For example, there are 14 partitions of
 3 with $(n+1)$
 copies of $n$:
\begin{align*}
&(3_3),\ (3_2),\ (3_1),\ (3_0),\  (2_2,1_1), \ (2_1,1_1),\ (
2_0,1_1),
\  (2_2,1_0), \  (2_1,1_0),\  (2_0,1_0), \\[6pt]
 &\qquad \qquad  (1_1,1_1,1_1),  \ (1_1,1_1,1_0),\
 (1_1,1_0,1_0), \ (1_0,1_0,1_0).
\end{align*}

\begin{figure}[htbp]
\centering
\includegraphics[width=15cm]{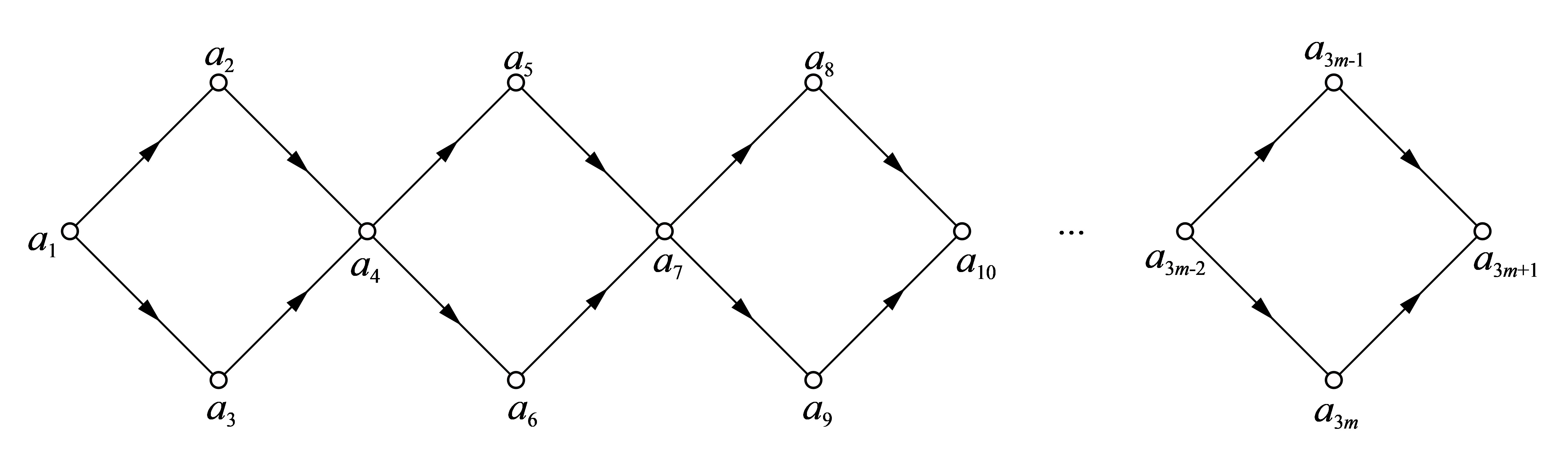}
\caption{A plane partition diamond of length $m$.}
\label{plane-partition}
\end{figure}

In 2001,   Andrews,   Paule and
 Riese \cite{Andrews-2001}
 introduced     plane
partition diamonds, which are partitions whose parts are
non-negative integers $a_i$ and $b_i
 $ which are placed at the nodes
of the graph given in the following graph (see Fig. \ref{plane-partition}),
with each directed edge
indicating $\geq$.

%
%
%
%
%
%
%
%
%
%
%
%
%

Additionally in \cite{Andrews-2022}, Andrews and Paule
 extended  Schmidt-type partitions to
partition diamonds. They only added up
  the summands at the links, namely,
   $a_1+a_4+a_7+\cdots$.
 Very recently, Andrews and Paule
 \cite{Andrews-2024}
  discussed   partition
diamonds with $(n +1)$
 copies of $n$ modified as follows.

 (I) The summands at the links (i.e.,
 $a_1, a_4, a_7, a_{10}, \ldots)
  $ are allowed to be
 $ (n +1) $ copies of $n$;

 (II) All other $a_i$ ($i\not\equiv 1
 \pmod 3)$ only
  have the subscript 0;

   (III) An arrow pointing from $a_i$
    to $a_j$
     means that $a_i\geq a_j$; the ``$\geq$"
      sign means that the weighted
       difference
    between parts is $\geq 0$, where
     the
     weighted difference between $m_i$
      and $n_j$ is defined by
      \[
((m_i-n_j)):=m-n-i-j.
      \]

\begin{figure}[h]
\centering
\includegraphics[width=15cm]{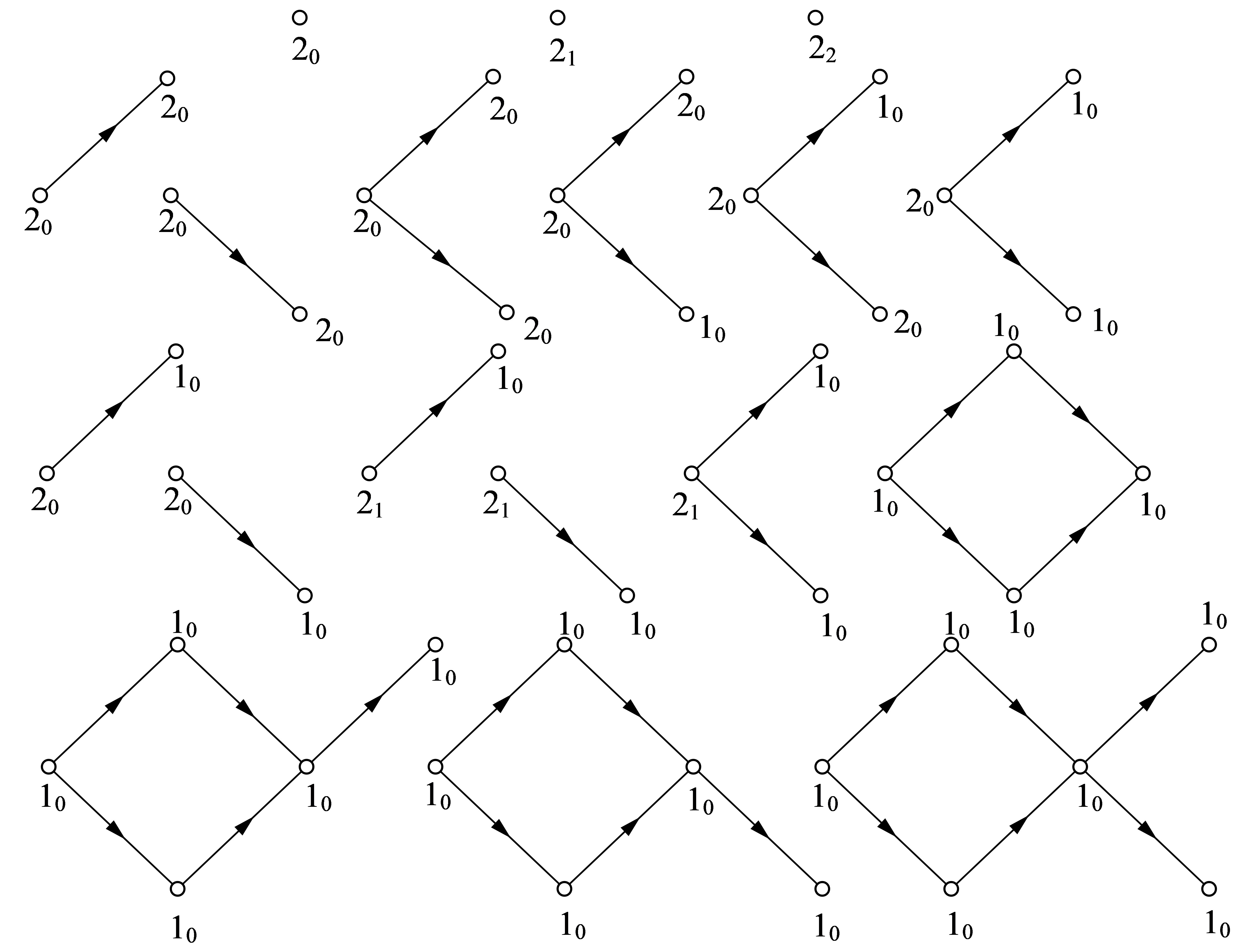}
\caption{Eighteen modified partition diamonds
   of 2 with $(n +1)$    copies of $n$.}
\label{2-partitions}
\end{figure}

  Let $PDN1 (N)$
 denote    the number of these
  modified partition diamonds with $(n +1)$
   copies of $n$  where summing the
   parts at the links gives $N$, namely,
    $a_1+a_4+\cdots=N$.
  For example, there are 18
   modified partition diamonds with
   $(n +1)$
   copies of $n$ where summing the
   parts at the links gives 2 (see Fig. \ref{2-partitions}).

  Andrews and Paule
  \cite{Andrews-2024} established
 the generating function for $PDN1(n)$  by means of partition
 analysis
 \begin{align}\label{1-1}
\sum_{n=0}^\infty PDN1(n) q^n=\frac{J_2^2
 }{ J_1^5 },
 \end{align}
where here and throughout this paper, we use the following
 $q$-series notation:
\begin{align*}
(a;q)_\infty &:=\prod_{n=0}^\infty(1-aq^n),
 \end{align*}
 and for any positive
  integer $k$,
  \[
J_k:=(q^k;q^k)_\infty.
  \]
Note that the generating
 function for $PDN1(n)$ is the
  reciprocal of the following
   identity which was independently discovered
by Ramanujan \cite{Ramanujan-1916}
 and Gordon \cite{Gordon}:
\[
\sum_{n=-\infty}^\infty 6(n+1)q^{3n(n+1)/2}
 =\frac{ J_1^5
 }{J_2^2}.
\]

Additionally, Andrews
       and Paule \cite{Andrews-2024}   proved
    various congruences  for $PDN1 (n)$
     by a combination of elementary methods:
     Radu's Ramanujan--Kolberg algorithm
      \cite{Radu} and Smoot's
      Mathematica package $\mathtt{RaduRK}$
      \cite{Smoot}. For example,
         they
   showed  that for $n\geq 0$,
 \begin{align}
  PDN1(25n+24)&\equiv  0 \pmod 5,
\label{1-2}\\[6pt]
PDN1(125n+74)&\equiv  0 \pmod {25},
\label{1-3}\\[6pt]
PDN1(125n+124)&\equiv  0 \pmod {25},
\label{1-4}\\[6pt]
   PDN1(7n+5)&\equiv  0 \pmod 7,\label{1-5}\\[6pt]
PDN1(49n+47)&\equiv  0 \pmod {49} \label{1-6}
\end{align}
and
\[
\sum_{j=0}^{5n+4}
 PDN1(25n+22-5j)PDN1(5j+1)\equiv 0
  \pmod 5.
\]

 In \cite{Lin},  Lin and Wang  considered
 another partition function
 which shares the same
  generating function with $PDN1(n)$.
 Lin and Wang \cite{Lin}
 presented human  proofs
  of  \eqref{1-2} and \eqref{1-5}.
 Congruences
   \eqref{1-3}, \eqref{1-4}
   and \eqref{1-6}
      were conjectured by Lin
      and Wang \cite{Lin}.
        Bian, Tang, Xia and Xue \cite{Bian}
         not only
          proved \eqref{1-3}
           and \eqref{1-4}, but also
        established
          some new congruences
        for $PDN1(n)$. For example, they
         proved
          that for
          $n\geq 0$,
\[
PDN1\left(\frac{23775n(3n+1)}{2}+974\right) \equiv
PDN1\left(\frac{23775n(3n+5)}{2}+974\right)
 \equiv 0 \pmod {25}.
\]
Very recently, Chen, Liu and Yao
  \cite{Chen} proved new congruences modulo $125$
  and $625$ for $PDN1(n)$. For example,
   they proved that for $n\geq 0$,
 \[
 PDN1(5^5 n+1849)
 \equiv  PDN1(5^5 n +3099)
\equiv  0 \pmod {5^3},
\]
and
\[
 PDN1(5^7n+46224)
\equiv PDN1(5^7 n+77474)
 \equiv 0 \pmod {5^4}.
 \]

At the end of their paper,
 Andrews and Paule \cite{Andrews-2024} stated
 ``Another task would be to determine
 infinite families of congruences
 similar to Ramanujan's classical $p(5^kn
 +d_k) \equiv 0
 \pmod {5^k}$ where $24d_k\equiv 1
  \pmod {5^{k}}$". Here, $p(n)$
 is the partition
  function which counts
   the number of partitions
    of $n$.
Motivated by Andrews and Paule's
   open problem, we shall prove
     congruences modulo arbitrary powers of
 $5$ and $7$ for $PDN1(n)$.
  In particular, we give an answer
   of Andrews and Paule's
   open problem.

The main results of this paper
 can be stated as follows.

\begin{theorem}\label{Th-1}
For any $n\geq 0$ and $\alpha\geq 1$,
\begin{align}\label{1-7}
PDN1\left(5^{2\alpha
 }n +\frac{23\cdot 5^{2\alpha}+1}{ 24}
 \right)\equiv 0 \pmod {5^{ \alpha }},
\end{align}
and
\begin{align}\label{1-8}
PDN1\left(5^{2\alpha+1}n +\frac{r\cdot 5^{2\alpha}+1}{ 24}
 \right)\equiv 0 \pmod {5^{ \alpha
 +1}},
\end{align}
where $r\in\{71,119\}$.
 \end{theorem}

\begin{theorem}\label{cong}
For any $n\geq 0$ and $\alpha\geq 1$,
\begin{align}
PDN1\left( 7^{2\alpha-1}n+\frac{17\cdot 7^{2\alpha-1}+1}{24}\right)
&\equiv 0\pmod{7^{\alpha}},\label{1-9}
\end{align}
and
\begin{align}
PDN1\left( 7^{2\alpha}n+\frac{23\cdot 7^{2\alpha}+1}{24}\right)
&\equiv 0\pmod{7^{\alpha+1}}.\label{1-10}
\end{align}
\end{theorem}

Congruences \eqref{1-7} and \eqref{1-8}
   are analogous
 to the
 following congruences modulo
  arbitrary powers of
 $5$ for $p(n)$:
\begin{align}\label{1-11}
p(5^{\alpha}n
 +\delta_{\alpha}) \equiv 0
 \pmod {5^{\alpha}},
  \end{align}
 where $\alpha$ is a
  positive
   integer and
     $\delta_{\alpha}$
     is   the reciprocal modulo
      $5^{\alpha}$ of 24. Congruence
       \eqref{1-11}
         was conjectured by Ramanujan
        \cite{Ramanujan-1} and was first proved by Watson
         \cite{Watson}. Hirschhorn and Hunt
         \cite{Hirschhorn-Hunt} gave a simple
         proof of \eqref{1-11}.  Watson
          \cite{Watson} also proved that
          for $\alpha\geq1$,
\begin{align}\label{1-12}
p(7^{2\alpha-1}n+\lambda_{2\alpha-1})
 \equiv 0 \pmod {7^{\alpha}}
\end{align}
and
\begin{align}\label{1-13}
p(7^{2\alpha }n+\lambda_{2\alpha })
 \equiv 0 \pmod {7^{\alpha+1}}.
\end{align}
Congruences \eqref{1-9} and \eqref{1-10}
  are analogous
 to Watson's  \eqref{1-12} and
  \eqref{1-13}.
Garvan \cite{Garvan}
  gave a simple proof of \eqref{1-12} and
  \eqref{1-13}.

In this paper, we give a
natural way to prove the partition congruences \eqref{1-7}--\eqref{1-10}
by increasing the rank of $\mathbb{Q}[t]$-modules,
which is the vector space over $\mathbb{Q}$
generated by the modular eta-quotients with poles only at infinity.
In the proof of the congruence \eqref{1-8},
the algorithm on the Ramanujan-type
identities presented by Chen, Du and Zhao \cite{CDZ19} is required.

The rest of this paper is organized as follows. In Sect.
\ref{Sect:Pre}, we collect some terminology and necessary results in
the theory of modular forms, and introduce how to find the modular
functions on the Ramanujan-type identities presented by Chen, Du and
Zhao \cite{CDZ19}. In Sect. \ref{sec-pf-th1}, we first establish the
generating functions of four infinite arithmetic progressions of
$PDN1(n)$ by utilizing the Atkin $U$-operator and the theory of
modular forms, and then give a proof of Theorem \ref{Th-1}. In Sect.
\ref{sect-pf-th2}, we investigate the congruences of $PDN1(n)$
modulo the powers of 7, and provide a proof of Theorem \ref{cong}.
In Appendix A, we shall give the expressions of some functions under
the action of the Atkin $U$-operator, which are used to derive the
generating functions in Sect.
 \ref{sec-pf-th1}.
The coefficients in the modular equation of the seventh order are
presented in Appendix B. Due to space constraints, we do not place
all of our computations in this article. Instead, we include a
Mathematica supplement, which can be found in
\url{https://github.com/tztgm1/PDN1-congruences-modulo-5}
and
\url{https://github.com/tztgm1/PDN1-congruences-modulo-7}.

\section{Preliminaries}\label{Sect:Pre}

In this section,
we collect some terminology and results in the theory of modular forms,
and state how to find the modular functions on
the Ramanujan-type identities of partition functions
defined in terms of eta-quotients
presented by Chen, Du and Zhao \cite{CDZ19}.

The full modular group is defined  by
\begin{align*}
\Gamma=\textrm{SL}_2(\mathbb{Z}): =\left\{\begin{pmatrix}a &b\\ c
&d\end{pmatrix}\colon
a,b,c,d\in\mathbb{Z},~\text{and}~ad-bc=1\right\},
\end{align*}
and for a positive integer $N$, the congruence subgroups are defined by
\begin{align*}
\Gamma_0(N)&:=\left\{\begin{pmatrix}a &b\\ c &d\end{pmatrix}\in\Gamma \colon c\equiv 0\pmod{N}\right\},\\[5pt]
\Gamma_1(N)&:=\left\{\begin{pmatrix}a & b\\ c &
d\end{pmatrix}\in\Gamma_0(N)\colon a\equiv d\equiv
1\pmod{N}\right\}.
\end{align*}
We denote by $\gamma$ the matrix $\begin{pmatrix}a &b\\ c &d \end{pmatrix}$, if not specified otherwise.
Let $\gamma$ act on $\tau\in \mathbb{C}$ by the linear fractional transformation
\begin{align*}
    \gamma\tau = \frac{a\tau+b}{c\tau+d} \qquad \text{and} \qquad
    \gamma\infty =\lim_{\tau\rightarrow \infty} \gamma\tau.
\end{align*}

Let $N$ and $k$ be positive integers and
$\mathbb{H}=\{\tau\in\mathbb{C}\colon \mathrm{Im}(\tau)>0\}$. A
holomorphic function $f\colon \mathbb{H}\rightarrow \mathbb{C}$ is
called a modular function of weight
 $k$ for $\Gamma'$
($\Gamma'=\Gamma_0(N) \ \text{or } \Gamma_1(N)$), if it satisfies
the following two conditions:

(1) for all $\gamma\in\Gamma'$, $f(\gamma\tau)=(c\tau+d)^kf(\tau)$;

(2) for any $\gamma\in\Gamma$, $(c\tau+d)^{-k}f(\gamma\tau)$ has a
Fourier expansion of the form
\begin{align*}
(c\tau+d)^{-k}f(\gamma\tau)=\sum_{n=n_\gamma}^\infty a(n)q_{w_{\gamma}}^n,
\end{align*}
where $a(n_\gamma)\neq 0$, $q_{w_{\gamma}}=e^{2\pi i\tau/w_{\gamma}}$,
and $w_\gamma$ is the minimal positive integer $h$ such that
\begin{align*}
\begin{pmatrix}1 &h\\ 0 &1\end{pmatrix}\in\gamma^{-1}\Gamma'\gamma.
\end{align*}
In particular, if $n_\gamma \geq 0$ for all $\gamma\in\Gamma$, then we call that $f$ is a modular form of weight $k$ for $\Gamma'$.

The modular function with weight 0 for $\Gamma'$ is referred to as a modular function for $\Gamma'$. For a modular function $f(\tau)$ of weight $k$ with respect to $\Gamma'$, the order of $f(\tau)$ at the cusp $a/c\in\mathbb{Q}\cup \{\infty\}$ is defined by
\begin{align*}
\mathrm{ord}_{a/c}(f)=n_\gamma
\end{align*}
for some $\gamma\in\Gamma$ such that $\gamma\infty =a/c$. It is known that $\mathrm{ord}_{a/c}(f)$ is well-defined, see
\cite[p. 72]{Diamond-Shurman-2005}.

Let $q=e^{2\pi i\tau}$ and $\tau\in\mathbb{H}$. The Dedekind eta-function $\eta(\tau)$ is defined by
\begin{align*}
\eta(\tau):=q^{\frac{1}{24}}\prod_{n=1}^\infty(1-q^n).
\end{align*}

The following   theorem \cite{GH93,Ligozat,Newman} is
  useful to verify whether an eta-quotient
is a modular form.

\begin{theorem}\label{GHN-THM}
  If
$f(\tau)=\prod_{\delta|N}\eta^{r_\delta}(\delta\tau)$ is an
eta-quotient with
$k=\frac{1}{2}\sum_{\delta|N}r_\delta\in\mathbb{Z}$, and satisfies
the following conditions:
\begin{align*}
\sum_{\delta|N}\delta r_\delta\equiv0\pmod{24}
\end{align*}
and
\begin{align*}
\sum_{\delta|N}\dfrac{N}{\delta}r_\delta\equiv0\pmod{24},
\end{align*}
then $f(\tau)$ fulfills
\begin{align*}
f(\gamma\tau)=\chi(d)(c\tau+d)^kf(\tau),
\end{align*}
for any $\gamma=\begin{pmatrix}a &b\\ c &d\end{pmatrix}\in\Gamma_0(N)$, where the character $\chi$ is defined
by $\chi(d)=\Big(\frac{(-1)^k s}{d}\Big)$, and $s=\prod_{\delta|N}\delta^{r_\delta}$.
\end{theorem}

The following theorem \cite{Biagioli,Ligozat,Martin}
 is the necessary criterion
  for determining orders of an eta-quotient
   at cusps.

\begin{theorem}\label{THM-cusp}
Let $c$, $d$ and $N$ be the positive integers with $d|N$ and $\gcd(c,d)=1$. If $f(\tau)$ is an
eta-quotient satisfying the conditions of Theorem \ref{GHN-THM}, then the order of
$f(\tau)$ at the cusp $c/d$ is
\begin{align*}
\dfrac{N}{24}\sum_{\delta|N}\dfrac{\gcd(d,\delta)^2r_\delta}{\gcd{\left(d,\frac{N}{d}\right)}d\delta}.
\end{align*}
\end{theorem}

For a positive integer $\delta$ and a residue class $g\pmod{\delta}$, the generalized Dedekind eta-function
$\eta_{\delta,g}(\tau)$ is defined as
\begin{align*}
\eta_{\delta,g}(\tau):=q^{\frac{\delta}{2}P_2{\left(\frac{g}{\delta}\right)}}
\prod_{\substack{n>0\\ n \equiv g\pmod{\delta}}}(1-q^n)
\prod_{\substack{n>0\\ n \equiv -g\pmod{\delta}}}(1-q^n),
\end{align*}
where
\begin{align*}
P_2(t):=\{t\}^2-\{t\}+\frac{1}{6}
\end{align*}
is the second Bernoulli function and $\{t\}$ is the fractional part
of $t$; see, for example  \cite{Robins-1994,Schoeneberg-1974}.
Notice that
\begin{align*}
\eta_{\delta, 0}(\tau)=\eta^2(\delta\tau)\qquad\text{and}\qquad\eta_{\delta, \frac{\delta}{2}}(\tau)
=\frac{\eta^2(\frac{\delta}{2}\tau)}{\eta^2(\delta\tau)}.
\end{align*}

Chen, Du and Zhao \cite{CDZ19}
devised an algorithm to compute the Ramanujan-type \mbox{identity} involving $a(mn+t)$ for any $m\geq 1$ and $0\leq t\leq m-1$
by finding suitable  modular functions for $\Gamma_1(N)$,
where $a(n)$ is defined by
\begin{align}
  \sum_{n=0}^{\infty}
a(n)q^n:=\prod_{\delta | M}
(q^\delta;q^\delta)^{r_\delta}_\infty,\label{radu-gf}
\end{align}
where $M$ is a positive integer and
$r_\delta$ are integers. Let
\begin{align*}
g_{m,t}(\tau)  = q^{\frac{t-\ell}{m}}\sum\limits_{n=0}^\infty a(mn+t)q^n,
\end{align*}
where
\[\ell=-\frac{1}{24}\sum\limits_{\delta|M}\delta r_\delta.\]

Let $\kappa=\gcd(m^2-1, 24)$.
Assume that $N$ satisfies the following conditions:
\begin{description}\label{con}
\setlength{\parskip}{2ex}
  \item[{1.}\label{con_N_1}] $M|N$.

  \item[{2.}\label{con_N_2}] $p| N$ for any prime $p| m$.

  \item[{3.}\label{con_N_6}] $\kappa N\sum\limits_{\delta|M}r_\delta \equiv 0 \pmod{8}$.

  \item[{4.}\label{con_N_7}] $\kappa mN^2\sum\limits_{\delta|M}\frac{r_\delta}{\delta}\equiv 0 \pmod{24}$.

  \item[{5.}\label{con_N_8}] $\frac{24m}{\gcd(\kappa \alpha(t), 24m)}\left| N\right.$,
  where $\alpha(t) = -\sum_{\delta|M}\delta r_{\delta}-24t$.

  \item[{6.}\label{con_N_9}] Let $\prod_{\delta|M}\delta^{|r_\delta|} = 2^zj$,
   where $z\in \mathbb{N}$ and $j$ is odd.
   If $2| m$, then $\kappa N\equiv 0 \pmod 4$ and $Nz\equiv 0 \pmod 8$,
    or $z\equiv 0 \pmod 2$ and $N(j-1)\equiv 0 \pmod 8$.

  \item[{7.}\label{con_N_10}]
  Let $\mathbb{S}_n = \{j^2 \pmod n \colon j\in\mathbb{Z}_n,\  \gcd(j, n)=1,\  j\equiv 1 \pmod{N}\}$.
  For any $s\in \mathbb{S}_{24mM}$,
 \[\frac{s-1}{24}\sum\limits_{\delta|M}\delta r_\delta + ts \equiv t \pmod m.\]
\end{description}

The following theorem \cite[Theorem 2.1]{CDZ19} gives
a necessary and sufficient condition for
the product of a generalized eta-quotient and
$g_{m,t}(\tau)$ to be a modular function for $\Gamma_1(N)$.

\begin{theorem}\label{con_F_modular_function}
For a given partition function $a(n)$ as defined by \eqref{radu-gf},
and for given integers $m$ and $t$, suppose that $N$ is a positive
integer satisfying the conditions 1--7. Let
\begin{align*}
  F(\tau)=\phi(\tau)\, g_{m,t}(\tau),
\end{align*}
where
\begin{align*}
\phi(\tau)=\prod_{\delta | N}\eta^{a_{\delta}}(\delta \tau)\,\prod_{{\delta|N\atop 0<g\leq \left\lfloor\delta/2\right\rfloor}}\eta_{\delta,g}^{a_{\delta,g}}(\tau),
\end{align*}
and $a_{\delta}$ and $a_{\delta,g}$ are integers.
Then $F(\tau)$ is a modular function with respect to $\Gamma_1(N)$ if and only if $a_{\delta}$ and $a_{\delta,g}$ satisfy the following conditions:

\begin{enumerate}
\setlength{\parskip}{2ex}
  \item[{\rm(1)}\label{F_con_1}] $\sum\limits_{\delta|N}a_\delta+\sum\limits_{\delta|M}r_\delta=0$,

  \item[{\rm(2)}\label{F_con_2}] $N \sum\limits_{\delta|N}\frac{a_\delta}{\delta}
                  +2N\sum\limits_{{\delta|N\atop 0<g\leq \left\lfloor{\delta}/{2}\right\rfloor}}\frac{a_{\delta,g}}{\delta}
                  +Nm\sum\limits_{\delta|M}\frac{r_\delta}{\delta}
                  \equiv0\pmod{24}$,

  \item[{\rm(3)}\label{F_con_3}]
  $\sum\limits_{\delta|N}\delta a_\delta
  +12 \sum\limits_{{\delta|N\atop 0<g\leq \left\lfloor{\delta}/{2}\right\rfloor}}\delta P_2\left(\frac{g}{\delta}\right){a_{\delta,g}}
  +m\sum\limits_{\delta|M}\delta r_\delta
  +\frac{(m^2-1)\alpha(t)}{m }
  \equiv0\pmod{24},$

  \item[{\rm(4)}\label{F_con_4}] For any integer $0<a<12N$ with $\gcd{(a,6)}=1$ and $a\equiv 1\pmod N$,
  $$\prod\limits_{\delta|N}\left(\frac{\delta}{a}\right)^{|a_\delta|}
  \prod\limits_{\delta|M}
  \left(\frac{m\delta}{a}\right)^{|r_\delta|} e^{\sum\limits_{\delta|N}\sum\limits_{g=1}^{{\tiny \left\lfloor\delta/2\right\rfloor}}\pi i\big(\frac{g}{\delta}-\frac{1}{2}\big)(a-1)a_{\delta,g}}=1.$$
\end{enumerate}
\end{theorem}

For any $\gamma\in\Gamma$ and $\lambda\in \mathbb{Z}_m$, let
\begin{align*}
p(\gamma, \lambda)  =  \frac{1}{24}\sum\limits_{\delta|M}\frac{\gcd^2(\delta(a+\kappa\lambda c), mc)}{\delta m}r_\delta,
\end{align*}
and for $\gamma\in\Gamma$, let
\begin{align}
p(\gamma) &= \min\{p(\gamma, \lambda)\colon \lambda=0, 1, \ldots, m-1\}.
\label{def_p_gamma-1}
\end{align}

The following theorem \cite[Theorem 2.7]{CDZ19} gives
estimates of the orders of $F(\tau)$ at cusps of $\Gamma_1(N)$.

\begin{theorem}\label{order_F-1}
For a given partition function $a(n)$ as defined by \eqref{radu-gf},
and for given integers $m$ and $t$,
let
\begin{align*}
  F(\tau)=\phi(\tau)\,g_{m,t}(\tau),
\end{align*}
where
\[\phi(\tau)=\prod_{\delta | N}\eta^{a_{\delta}}(\delta \tau)\prod_{{\delta|N\atop 0<g\leq \left\lfloor\delta/2\right\rfloor}}\eta_{\delta,g}^{a_{\delta,g}}(\tau),\] $a_{\delta}$ and $a_{\delta,g}$ are integers.
Assume that $F(\tau)$ is a modular function for $\Gamma_1(N)$.
Let $\{s_1, s_2,\ldots,s_\epsilon\}$ be a complete set of inequivalent cusps of $\Gamma_1(N)$,
and for each $1\leq i \leq \epsilon$,
let $\gamma_i\in\Gamma$ be such that $\gamma_i\infty = s_i$.
Then
\begin{align*}
\mathrm{ord}_{s_i}(F(\tau))\ge w_{\gamma_i}\,(p(\gamma_i)+p^*(\gamma_i)),
\end{align*}
where $p(\gamma)$ is given by \eqref{def_p_gamma-1} and $p^*(\gamma)$ is defined as
\begin{align*}
p^*(\gamma)&=\frac{1}{24}\sum_{\delta|N}\frac{\gcd^2{(\delta,c)}}{\delta}a_\delta
+\frac{1}{2}\sum_{{\delta|N\atop 0<g\leq \left\lfloor\delta/2\right\rfloor}}\frac{\gcd^2{(\delta,c)}}{\delta}P_2\left(\frac{a g}{\gcd{(\delta,c)}}\right) a_{\delta,g}.
\end{align*}
\end{theorem}

The following theorem of Sturm \cite[Theorem~1]{Sturm-1987} plays
an important role in proving congruences using the theory of modular forms.

\begin{theorem}\label{Sturm-theorem}
Let $\Gamma'$ be a congruence subgroup of $\Gamma$, and let $k$ be
an integer and $g(\tau)=\sum\limits_{n=0}^\infty c(n)q^n$ be  a
modular form of weight $k$ for $\Gamma'$. For any given positive
integer $u$, if $c(n)\equiv0\pmod{u}$ holds for all $n\leq
\frac{k}{12}[\Gamma:\Gamma']$, then $c(n)\equiv0\pmod{u}$ holds for
any $n\geq 0$.
\end{theorem}

There is an explicit formula for the index \cite[p.~13]{Diamond-Shurman-2005}:
\begin{align*}
[\Gamma: \Gamma_1(N)]=N^2\cdot\prod_{\substack{p|N\\p~\text{is~prime}}}{\left(1-\frac{1}{p^2}\right)}.
\end{align*}

For a function $f\colon \mathbb{H}\rightarrow \mathbb{C}$ and a prime $p$, the Atkin $U$-operator $U_p$ is
defined as
\begin{align*}
U_p(f)(\tau): =\frac{1}{p}\sum_{\lambda=0}^{p-1}
f{\left(\frac{\tau+\lambda}{p}\right)},
\end{align*}
where $\tau\in\mathbb{H}$.
Suppose $f$ has an expansion $f=\sum_{n=-\infty}^\infty a(n)q^n$.
It is easy to see that
\begin{align*}
U_p{\left(\sum_{n=-\infty}^\infty a(n)q^n\right)}=\sum_{n=-\infty}^\infty a(pn)q^n,
\end{align*}
where $q=e^{2\pi i\tau}$. It is known that if $f$ is a modular function for $\Gamma_0(pN)$, where
$p|N$, then $U_p(f)$ is a modular function for $\Gamma_0(N)$. Gordon and Hughes \cite[Theorem 4]{GH81} gave the lower
bounds on the order of $U_p(f)$ at cusps of $\Gamma_0(N)$ in terms of that of $f$ at cusps of $\Gamma_0(pN)$. Denote
by $\pi(n)$ the $p$-adic order of $n$.

\begin{lemma}\label{GH-ord}
Suppose that $f$ is a modular function for $\Gamma_0(pN)$, where $p|N$. Let $r=a/c$ be a cusp of $\Gamma_0(N)$,
where $c|N$ and $\gcd(a,c)=1$. Then
\begin{align*}
\mathrm{ord}_r\big(U_p(f)\big)\geq
\begin{cases}
\frac{1}{p}\mathrm{ord}_{r/p}(f), & \text{if } \pi(c)\geq \frac{1}{2}\pi(N);\\[5pt]
\mathrm{ord}_{r/p}(f), & \text{if } 0<\pi(c)<\frac{1}{2}\pi(N);\\[5pt]
\min_{0\leq \lambda\leq p-1}\mathrm{ord}_{(r+\lambda)/p}(f), & \text{if } \pi(c)=0.
\end{cases}
\end{align*}
\end{lemma}

\section{Congruences modulo the powers of $5$}\label{sec-pf-th1}

\subsection{The generating functions}\label{Sect-gf}

In this subsection, we shall build a modular equation of the fifth order, and then derive the expressions of the generating functions of
four infinite arithmetic progressions on $PDN1(n)$.

Recall the generating function of $PDN1(n)$
\begin{align*}
\sum_{n=0}^\infty PDN1(n) q^n=\frac{J_2^2}{J_1^5}.
 \end{align*}
Consider the function
\begin{align*}
A:=\frac{\eta^2(2 \tau ) \eta^5(25\tau
)}{\eta^5(\tau)\eta^2(50\tau)}.
\end{align*}
We recursively define the $U$-sequence $(L_\alpha)_{\alpha\geq 0}$ by $L_0=1$, and for any $\alpha\geq 1$,
\begin{align*}
L_{2\alpha-1} &=U_5(A\!\,L_{2\alpha-2}),\\
L_{2\alpha} &=U_5(L_{2\alpha-1}).
\end{align*}
Then for any $\alpha\geq 1$,
\begin{align}
L_{2\alpha-1}&=\frac{J_5^5}{J_{10}^2}q\sum_{n=0}^\infty
PDN1\left( 5^{2\alpha-1}n+\frac{95\cdot 5^{2\alpha-2}+1}{24}\right)q^n,\label{L-odd-gf}\\[5pt]
L_{2\alpha}&=
\frac{J_1^5}{J_{2}^2}q\sum_{n=0}^\infty PDN1
\left(
5^{2\alpha}n+\frac{23\cdot 5^{2\alpha}+1}{24}\right)q^n.\label{L-even-gf}
\end{align}
It follows  from Theorem \ref{GHN-THM} that $A$ is a modular
function for $\Gamma_0(50)$.
Moreover, by Theorem \ref{THM-cusp}, we obtain the orders of $A$
at the cusps of $\Gamma_0(50)$:
\begin{align*}
&\mathrm{ord}_0(A)=-8,\  \mathrm{ord}_{1/2}(A)=-1,\
\mathrm{ord}_{1/25}(A)=8,\  \mathrm{ord}_{\infty}(A)=1,\\[5pt]
&\mathrm{ord}_{1/5}(A)=\mathrm{ord}_{2/5}(A)
=\mathrm{ord}_{3/5}(A)=\mathrm{ord}_{4/5}(A)=0,\\[5pt]
&\mathrm{ord}_{1/10}(A)=\mathrm{ord}_{3/10}(A)
=\mathrm{ord}_{7/10}(A)=\mathrm{ord}_{9/10}(A)=0.
\end{align*}
Thus, $U_5(A)$ is a modular function for $\Gamma_0(10)$,
and, by Lemma \ref{GH-ord}, we have
\begin{align}
\mathrm{ord}_0(U_5(A)) \geq -8,\quad \mathrm{ord}_{1/2}(U_5(A))\geq -1, \quad \mathrm{ord}_{1/5}(U_5(A)) \geq 2,
\quad \mathrm{ord}_\infty(U_5(A)) \geq 1.\label{order-U-A}
\end{align}
Hence, $L_\alpha (\alpha\geq 1)$ are modular functions
for $\Gamma_0(10)$.
In order to obtain the expressions of $L_\alpha$, we define
\begin{align*}
t:=\frac{\eta^3(\tau)
\eta(5\tau)}{\eta(2\tau)\eta^3(10\tau)},
\quad z:=\frac{\eta(\tau)\eta^3(2\tau)\eta^3(5\tau)}{\eta^{7}(10\tau)},
\quad u:=z^{-1}=
\frac{\eta^{7}(10\tau)}{\eta(\tau)\eta^3(2\tau)\eta^3(5\tau)}.
\end{align*}
It follows from Theorem \ref{GHN-THM} that
$t, z$ and $u$ are modular functions for $\Gamma_0(10)$.
By Theorem \ref{THM-cusp}, we obtain the orders of $t$ and $z$
at the cusps of $\Gamma_0(10)$:
\begin{align}
\mathrm{ord}_0(t) &=1, &\mathrm{ord}_{1/2}(t) &=0, &\mathrm{ord}_{1/5}(t) &=0, &\mathrm{ord}_\infty(t) &=-1,\label{order-t}\\
\mathrm{ord}_0(z) &=1, &\mathrm{ord}_{1/2}(z) &=1, &\mathrm{ord}_{1/5}(z) &=0, &\mathrm{ord}_\infty(z) &=-2.\label{order-z}
\end{align}
So the function $t$ is a hauptmodul of modular functions on
$\Gamma_0(10)$.

Denote by $E^\infty(N)$ the set of modular eta-quotients with poles only at infinity for $\Gamma_0(N)$, and denote by
$\langle E^\infty(N) \rangle_\mathbb{Q}$ the vector space over $\mathbb{Q}$ generated by the elements of $E^\infty(N)$.
Therefore, 
\begin{align*}
\langle E^\infty(10) \rangle_\mathbb{Q}=\mathbb{Q}[t].
\end{align*}
Using the algorithm for proving the
theta-function identities developed by Frye and Garvan \cite{FG} or
the zero recognition of modular functions presented by Paule and
Radu \cite{PR21}, we derive
\begin{align}
t^2=z - 5 t=u^{-1} - 5 t.\label{t-2}
\end{align}
So, $\langle E^\infty(10) \rangle_\mathbb{Q}$ can be rewritten as
\begin{align}
\langle E^\infty(10) \rangle_\mathbb{Q}=\mathbb{Q}[z]+\mathbb{Q}[z]t.
\label{E-10-z-t}
\end{align}
In order to annihilate the poles of $U_5(A)$ at the other cusps
but $\infty$, we can allow $U_5(A)$ to multiply by $z^8$
such that $z^8U_5(A)$ has a pole only at $\infty$.
To obtain a useful expression of $U_5(A)$,
we use the expression \eqref{E-10-z-t} of $\langle E^\infty(10) \rangle_\mathbb{Q}$.

The following theorem gives the expressions of $L_{\alpha}$
for $\alpha\geq 1$.

\begin{theorem}\label{L-expression}
For any $\beta\geq 0$, we have
\begin{align}
L_{4\beta+1}&=5^{2\beta}a_{4\beta+1,0}(3u+ut+u^2t)
+5^{2\beta+2}a_{4\beta+1,1}u+5^{2\beta+2}a_{4\beta+1,2}u^2+
5^{2\beta+1}\sum_{i=3}^{(5^{4\beta+2}-1)/3} 5^{\lfloor
\frac{5i}{6}\rfloor} a_{4\beta+1,i}u^i\notag\\[5pt]
&\quad+5^{2\beta+2}b_{4\beta+1,1}ut+5^{2\beta+2}b_{4\beta+1,2}u^2t+
5^{2\beta+1}\sum_{i=3}^{(5^{4\beta+2}-1)/3} 5^{\lfloor
\frac{5i}{6}\rfloor} b_{4\beta+1,i}u^it,\label{L-4n-1}\\[5pt]
L_{4\beta+2}&=5^{2\beta+1}a_{4\beta+2,0}(ut-u^2t)
+5^{2\beta+2}a_{4\beta+2,1}u+5^{2\beta+2}a_{4\beta+2,2}u^2+
5^{2\beta+1}\sum_{i=3}^{(5^{4\beta+3}-5)/3} 5^{\lfloor
\frac{5i}{6}\rfloor} a_{4\beta+2,i}u^i\notag\\[5pt]
&\quad+5^{2\beta+2}b_{4\beta+2,1}ut+5^{2\beta+2}b_{4\beta+2,2}u^2t+
5^{2\beta+1}\sum_{i=3}^{(5^{4\beta+3}-5)/3} 5^{\lfloor
\frac{5i}{6}\rfloor} b_{4\beta+2,i}u^it,\label{L-4n-2}\\[5pt]
L_{4\beta+3}&=5^{2\beta+1}a_{4\beta+3,0}(7u-ut+4u^2t-15u^2)
+5^{2\beta+3}a_{4\beta+3,1}u+5^{2\beta+3}a_{4\beta+3,2}u^2\notag\\[5pt]
&\quad+ 5^{2\beta+2}\sum_{i=3}^{(5^{4\beta+4}-1)/3} 5^{\lfloor
\frac{5i}{6}\rfloor} a_{4\beta+3,i}u^i
+5^{2\beta+3}b_{4\beta+3,1}ut+5^{2\beta+3}b_{4\beta+3,2}u^2t
\notag\\[5pt]
&\quad+5^{2\beta+3}b_{4\beta+3,3}u^3t + 5^{2\beta+2}\sum_{i=4}^{
(5^{4\beta+4}-1)/3} 5^{\lfloor \frac{5i}{6}\rfloor}
b_{4\beta+3,i}u^it,
\label{L-4n-3}\\[5pt]
L_{4\beta+4}&=5^{2\beta+2}a_{4\beta+4,0}(u-u^2t)
+5^{2\beta+3}a_{4\beta+4,1}u+5^{2\beta+3}a_{4\beta+4,2}u^2+
5^{2\beta+2}\sum_{i=3}^{(5^{4\beta+5}-5)/3} 5^{\lfloor
\frac{5i}{6}\rfloor} a_{4\beta+4,i}u^i\notag\\[5pt]
&\quad+5^{2\beta+3}b_{4\beta+4,1}ut+5^{2\beta+3}b_{4\beta+4,2}u^2t+
5^{2\beta+2}\sum_{i=3}^{(5^{4\beta+5}-5)/3} 5^{\lfloor
\frac{5i}{6}\rfloor} b_{4\beta+4,i}u^it,\label{L-4n-4}
\end{align}
where $a_{i,j}$ and $b_{i,j}$ are integers.
\end{theorem}


To prove  Theorem \ref{L-expression}, we require the following
fifth order modular equation.

\begin{lemma}\label{modular-equation}
Let
\begin{align}
\sigma_1(q)&=1125 u + 30 t u + 114625 u^2 + 10625 t u^2 + 3137500 u^3 +
 435000 t u^3 + 31015625 u^4\notag\\[5pt]
 &\quad + 5359375 t u^4 + 100000000 u^5 +
 19921875 t u^5,\notag\\[5pt]
\sigma_2(q)&=635 u + 5 t u + 96500 u^2 + 8825 t u^2 + 2925000 u^3 + 406875 t u^3 +
 30250000 u^4\notag\\[5pt]
 &\quad + 5250000 t u^4 + 100000000 u^5 + 20000000 t u^5,\notag\\[5pt]
\sigma_3(q)&=160 u + 31125 u^2 + 2700 t u^2 + 1025625 u^3 + 140625 t u^3 +
 11062500 u^4 + 1912500 t u^4 \notag\\[5pt]
 &\quad+ 37500000 u^5 + 7500000 t u^5,\notag\\[5pt]
\sigma_4(q)&=20 u + 4525 u^2 + 375 t u^2 + 160000 u^3 + 21625 t u^3 +
 1796875 u^4 + 309375 t u^4  \notag\\[5pt]
 &\quad+ 6250000 u^5+ 1250000 t u^5,\notag\\[5pt]
\sigma_5(q)&=u + 250 u^2 + 20 t u^2 + 9375 u^3 + 1250 t u^3 + 109375 u^4 +
 18750 t u^4 + 390625 u^5 + 78125 t u^5.\label{sigma-5}
\end{align}
Then
\begin{align}
u^5=\sigma_1(q^5) u^4+\sigma_2(q^5) u^3+\sigma_3(q^5) u^2
+\sigma_4(q^5) u +\sigma_5(q^5).\label{modular-equation-eq}
\end{align}
\end{lemma}

\pf
Since $u$ is a modular function for $\Gamma_0(10)$, as well as for
$\Gamma_0(50)$, the functions on both sides of
\eqref{modular-equation-eq} are modular functions with respect to
$\Gamma_0(50)$. So, \eqref{modular-equation-eq} can be shown by
using the algorithm of proving the
theta-function identities of Frye and Garvan \cite{FG} or
the zero recognition of modular functions of Paule and
Radu \cite{PR21}.
\qed

From Lemma  \ref{modular-equation}, we have

\begin{corollary}
For any function $f\colon \mathbb{H}\rightarrow \mathbb{C}$
and any $n\geq 5$, we have
\begin{align}
U_5(fu^n)=\sum_{i=1}^5 \sigma_i(q) U_5(fu^{n-i}).\label{u-f-u}
\end{align}
\end{corollary}

Utilizing the modular equation \eqref{modular-equation-eq}, we deduce
  the following lemma.

\begin{lemma}\label{U-A-u-n}
For any $n\geq 1$, let $v_0(n)=\left\lceil
\frac{n+1}{5}\right\rceil, v_1(n)=\left\lceil
\frac{n+3}{5}\right\rceil, \mu(n,i)=\left\lfloor
\frac{5i-n+1}{6}\right\rfloor$. Then there exist integers
$c_{n,i}$ and $d_{n,i}$ such that
\begin{align}
U_5(Au^n)=
\begin{cases}
&c_{n,v_0(n)} u^{v_0(n)}+\sum\limits_{i=v_0(n)+1}^{5n+8}
5^{\mu(n,i)} c_{n,i}u^i
+d_{n,v_1(n)} u^{v_1(n)} t +d_{n,v_1(n)+1} u^{v_1(n)+1} t\\[5pt]
&\quad +\sum\limits_{i=v_1(n)+2}^{5n+8}
5^{\mu(n,i)} d_{n,i}u^it,
\qquad\qquad \qquad\qquad \ \qquad\text{if } n \equiv 0,1,2 \pmod{5},\\[5pt]
&c_{n,v_0(n)} u^{v_0(n)}+c_{n,v_0(n)+1} u^{v_0(n)+1}
+\sum\limits_{i=v_0(n)+2}^{5n+8} 5^{\mu(n,i)} c_{n,i}u^i
+d_{n,v_1(n)} u^{v_1(n)} t \\[5pt]
&\quad +5 d_{n,v_1(n)+1} u^{v_1(n)+1} t
+\sum\limits_{i=v_1(n)+2}^{5n+8} 5^{\mu(n,i)} d_{n,i}u^it, \qquad
\text{if } n \equiv 3,4 \pmod{5}.
\end{cases}\label{U-A-u}
\end{align}
\end{lemma}

\pf
We prove \eqref{U-A-u}  by induction on $n$. It is easy
to check  that \eqref{U-A-u} is true for $1\leq n \leq 5$
in our Mathematica supplement.
Suppose that \eqref{U-A-u} holds for all
 integers less than $n$ with $n\geq6$.
Next we prove that it holds for $n$.
According to the residues of
$n$ modulo 5, there are five cases to be considered:
$n\equiv 0, 1, 2, 3, 4 \pmod{5}$.
Here, we only prove the case $n\equiv 0\pmod{5}$, the rest
can be shown analogously. Let $n=5 k$, where $k\geq 2$.
It follows from \eqref{u-f-u} that
\begin{align*}
U_5(Au^n)&=\sum_{i=1}^5 \sigma_i(q) U_5(Au^{n-i}).
\end{align*}
We prove that for each $1\leq i\leq 5$,
$\sigma_i(q) U_5(Au^{n-i})$ can be expressed as the form of
\begin{align}
a_{k+1}u^{k+1}
+\sum_{j=k+2}^{25k+8}5^{\lfloor\frac{5j-5k+1}{6}\rfloor}a_j
u^j +b_{k+1}u^{k+1}t+b_{k+2}u^{k+2}t
+\sum_{j=k+3}^{25k+8}5^{\lfloor\frac{5j-5k+1}{6}\rfloor}b_j
u^jt,\label{U-e-form}
\end{align}
where $a_j, b_j$ are some integers.
Next, we only prove that it
holds for $\sigma_5(q) U_5(Au^{n-5})$, and the other terms can be shown in the same manner.

By means of \eqref{sigma-5}, let
\[\sigma_5(q)=A_0+A_1t,\]
where
\begin{align}
A_0& :=u + 250 u^2  + 9375 u^3  + 109375 u^4  + 390625 u^5 ,\label{A-0}
 \\[6pt]
 A_1& :=20 u^2+ 1250 u^3+ 18750 u^4+ 78125 u^5.\label{A-1}
\end{align}
By the inductive hypothesis, we have
\begin{align*}
\sigma_5(q) U_5(Au^{n-5})&=
\sigma_5(q)\bigg(a_{5k-5,k} u^{k}
+\sum\limits_{i=k+1}^{25k-17} 5^{\lfloor \frac{5i-5k+6}{6}\rfloor}
a_{5k-5,i} u^i+b_{5k-5,k}u^kt+b_{5k-5,k+1}u^{k+1}t\\[5pt]
&\quad + \sum\limits_{i=k+2}^{25k-17} 5^{\lfloor \frac{5i-5k+6}{6}\rfloor}b_{5k-5,i} u^it\bigg)\\[5pt]
&=A_0\bigg(a_{5k-5,k} u^{k}
+\sum\limits_{i=k+1}^{25k-17} 5^{\lfloor \frac{5i-5k+6}{6}\rfloor}
a_{5k-5,i} u^i\bigg) \\[5pt]
&\quad+\bigg(A_0 \bigg(
b_{5k-5,k}u^k+b_{5k-5,k+1}u^{k+1}+ \sum\limits_{i=k+2}^{25k-17} 5^{\lfloor \frac{5i-5k+6}{6}\rfloor}b_{5k-5,i} u^i
\bigg)\nonumber\\[5pt]
&\quad
  +A_1 \bigg(a_{5k-5,k} u^{k}
+\sum\limits_{i=k+1}^{25k-17} 5^{\lfloor \frac{5i-5k+6}{6}\rfloor}
a_{5k-5,i} u^i\bigg)\bigg)t   \nonumber\\[5pt]
&\quad +A_1 \bigg(
b_{5k-5,k}u^k+b_{5k-5,k+1}u^{k+1}+ \sum\limits_{i=k+2}^{25k-17} 5^{\lfloor \frac{5i-5k+6}{6}\rfloor}b_{5k-5,i} u^i
\bigg)t^2.
\end{align*}
After substituting \eqref{A-0} and \eqref{A-1}
into the above identity, we obtain that $\sigma_5(q) U_5(Au^{n-5})$
can be expressed as
\begin{align*}
&a_{5k-5,k}u^{k+1}
+\sum_{i=k+2}^{25k-12} 5^{\lfloor\frac{5i-5k+1}{6} \rfloor}a_iu^i
+b_{k+1}u^{k+1}t+b_{k+2}u^{k+2}t+\sum_{i=k+3}^{25k-12} 5^{\lfloor\frac{5i-5k+1}{6} \rfloor}b_iu^it\\[5pt]
&\quad
+5c_{k+2}u^{k+2}t^2+5c_{k+3}u^{k+3}t^2+\sum_{i=k+4}^{25k-12} 5^{\lfloor\frac{5i-5k+2}{6} \rfloor}c_iu^{i+2}t^2,
\end{align*}
where $a_i, b_i, c_i$ are integers.
From \eqref{t-2}, we obtain that there exist
integers $a_i', b_i'$ such that
\begin{align*}
&\quad \sigma_5(q) U_5(Au^{n-5})\\[5pt]
&=
a_{k+1}'u^{k+1}
+\sum_{i=k+2}^{25k-12} 5^{\lfloor\frac{5i-5k+1}{6} \rfloor}a'_iu^i
+b'_{k+1}u^{k+1}t+b'_{k+2}u^{k+2}t+\sum_{i=k+3}^{25k-12} 5^{\lfloor\frac{5i-5k+1}{6} \rfloor}b'_iu^it.
\end{align*}
Therefore, \eqref{U-e-form} holds for $\sigma_5(q) U_5(Au^{n-5})$.
Hence, \eqref{U-A-u} is true for $n\equiv 0\pmod{5}$.
\qed

Analogously, we can derive the expressions of $U_5(Au^nt), U_5(u^n), U_5(u^nt)$.

\begin{lemma}\label{U-A-u-n-t}
For any $n\geq 1$, let $v_0(n)=\left\lceil \frac{n}{5}\right\rceil,
v_1(n)=\left\lceil \frac{n+3}{5}\right\rceil, \mu(n,i)=\left\lfloor
\frac{5i-n+1}{6}\right\rfloor$. Then there exist integers
$c_{n,i}$ and $d_{n,i}$ such that
\begin{align*}
U_5(Au^nt)=
\begin{cases}
&c_{n,v_0(n)} u^{v_0(n)}+\sum\limits_{i=v_0(n)+1}^{5n+3}
5^{\mu(n,i)} c_{n,i}u^i
+d_{n,v_1(n)} u^{v_1(n)} t +d_{n,v_1(n)+1} u^{v_1(n)+1} t\\[5pt]
&\quad +\sum\limits_{i=v_1(n)+2}^{5n+3}
5^{\mu(n,i)} d_{n,i}u^it,
\qquad \qquad \qquad \qquad \qquad \qquad \quad   \text{if } n \equiv 1,2 \pmod{5},\\[5pt]
&c_{n,v_0(n)} u^{v_0(n)}+c_{n,v_0(n)+1} u^{v_0(n)+1}
+\sum\limits_{i=v_0(n)+2}^{5n+3}
5^{\mu(n,i)} c_{n,i}u^i
+d_{n,v_1(n)} u^{v_1(n)} t \\[5pt]
&\quad +5 d_{n,v_1(n)+1} u^{v_1(n)+1} t
+\sum\limits_{i=v_1(n)+2}^{5n+3}
5^{\mu(n,i)} d_{n,i}u^it,
\qquad \qquad \text{if } n \equiv 3,4 \pmod{5},\\[5pt]
&c_{n,v_0(n)} u^{v_0(n)}+c_{n,v_0(n)+1} u^{v_0(n)+1}
+\sum\limits_{i=v_0(n)+2}^{5n+3}
5^{\mu(n,i)} c_{n,i}u^i
+d_{n,v_1(n)} u^{v_1(n)} t \\[5pt]
&\quad +\sum\limits_{i=v_1(n)+1}^{5n+3}
5^{\mu(n,i)} d_{n,i}u^it,
\qquad \qquad \qquad \qquad \qquad \qquad \qquad
\text{if } n \equiv 0 \pmod{5}.
\end{cases}
\end{align*}
\end{lemma}

\begin{lemma}\label{U-u-n}
For any $n\geq 1$,
let $v_0(n)=\left\lceil \frac{n}{5}\right\rceil,
v_1(n)=\left\lceil \frac{n+3}{5}\right\rceil,
\mu(n,i)=\left\lfloor \frac{5i-n+2}{6}\right\rfloor$.
Then there exist integers $c_{n,i}$ and $d_{n,i}$ such that
\begin{align*}
U_5(u^n)=
\begin{cases}
&\sum\limits_{i=v_0(n)}^{5n} 5^{\mu(n,i)}c_{n,i} u^i +d_{n,v_1(n)}
u^{v_1(n)}t +\sum\limits_{i=v_1(n)+1}^{5n} 5^{\mu(n,i)}d_{n,i} u^it,
\qquad  \ \,  \text{if }
n\equiv 1\pmod{5},\\
&c_{n,v_0(n)} u^{v_0(n)}
+\sum\limits_{i=v_0(n)+1}^{5n}
5^{\mu(n,i)}c_{n,i} u^i
\\
&\quad +d_{n,v_1(n)} u^{v_1(n)}t +\sum\limits_{i=v_1(n)+1}^{5n}
 5^{\mu(n,i)}d_{n,i} u^it,
\qquad\qquad\qquad\qquad
\qquad \text{if } n\equiv 2\pmod{5},\\
&c_{n,v_0(n)} u^{v_0(n)} +\sum\limits_{i=v_0(n)+1}^{5n}
5^{\mu(n,i)}c_{n,i} u^i + \sum\limits_{i=v_1(n)}^{5n}
5^{\mu(n,i)}d_{n,i} u^it, \qquad \  \ \, \text{if } n\equiv 0, 3,
4\pmod{5}.
\end{cases}
\end{align*}
\end{lemma}

\begin{lemma}\label{U-u-n-t}
For any $n\geq 1$,
let $v_0(n)=\left\lceil \frac{n}{5}\right\rceil,
v_1(n)=\left\lceil \frac{n+2}{5}\right\rceil,
\mu(n,i)=\left\lfloor \frac{5i-n+3}{6}\right\rfloor$.
Then there exist integers $c_{n,i}$ and $d_{n,i}$ such that
\begin{align*}
U_5(u^nt)=
\begin{cases}
&\sum\limits_{i=v_0(n)}^{5n-1} 5^{\mu(n,i)}c_{n,i} u^i
+d_{n,v_1(n)} u^{v_1(n)}t
+\sum\limits_{i=v_1(n)+1}^{5n} 5^{\mu(n,i)}d_{n,i} u^it,
\qquad \text{if } n\equiv 1,2,3\pmod{5},\\
&c_{n,v_0(n)} u^{v_0(n)}
+\sum\limits_{i=v_0(n)+1}^{5n-1} 5^{\mu(n,i)}c_{n,i} u^i+ \sum\limits_{i=v_1(n)}^{5n} 5^{\mu(n,i)}d_{n,i} u^it,
 \qquad \text{if } n\equiv 4\pmod{5},\\
&c_{n,v_0(n)} u^{v_0(n)}
+\sum\limits_{i=v_0(n)+1}^{5n-1} 5^{\mu(n,i)}c_{n,i} u^i\\
&\quad +d_{n,v_1(n)} u^{v_1(n)}t
+\sum\limits_{i=v_1(n)+1}^{5n} 5^{\mu(n,i)}d_{n,i} u^it,
 \qquad \qquad \qquad \qquad \quad \text{if } n\equiv 0\pmod{5}.
\end{cases}
\end{align*}
\end{lemma}

Now, we are ready to prove  Theorem \ref{L-expression}.

\noindent{\it Proof  of Theorem \ref{L-expression}}. We prove the
statement by induction on $\beta$. In view of  \eqref{U-A},
\begin{align*}
L_1&=U_5(A)\\[5pt]
&=-(3 u + u t + u^2 t)+
8675 u + 3277600 u^2 + 295664000 u^3 + 10852800000 u^4\\[5pt]
 &\quad +
 195776000000 u^5 + 1840000000000 u^6 + 8640000000000 u^7 +
 16000000000000 u^8\\[5pt]
 &\quad+
 150 ut + 193825 u^2t + 28016000 u^3t + 1343680000 u^4t +
 28940800000 u^5t \\[5pt]
 &\quad + 309120000000 u^6t+ 1600000000000 u^7t +
 3200000000000 u^8t,
\end{align*}
which implies that \eqref{L-4n-1} holds for $\beta=0$.

We assume \eqref{L-4n-1} holds for $\beta\geq 0$.
Next, we prove \eqref{L-4n-2}--\eqref{L-4n-4} hold for $\beta$,
and \eqref{L-4n-1} holds for $\beta+1$. Applying the
$U_5$-operator to $L_{4\beta+1}$ yields
\begin{align*}
L_{4\beta+2}&=U_5(L_{4\beta+1})\\[5pt]
&=5^{2\beta}a_{4\beta+1,0}U_5(3u+ut+u^2t)
+5^{2\beta+2}a_{4\beta+1,1}U_5(u)+5^{2\beta+2}a_{4\beta+1,2}U_5(u^2)\\[5pt]
&\quad+
5^{2\beta+1}\sum_{i=3}^{(5^{4\beta+2}-1)/3} 5^{\lfloor
\frac{5i}{6}\rfloor} a_{4\beta+1,i}U_5(u^i)+5^{2\beta+2}b_{4\beta+1,1}U_5(ut)
+5^{2\beta+2}b_{4\beta+1,2}U_5(u^2t)\\[5pt]
&\quad+
5^{2\beta+1}\sum_{i=3}^{(5^{4\beta+2}-1)/3} 5^{\lfloor
\frac{5i}{6}\rfloor} b_{4\beta+1,i}U_5(u^it).
\end{align*}
In views of \eqref{U-3-u} and Lemmas \ref{U-u-n} and \ref{U-u-n-t},
there exist integers $a_{4\beta+2,i}$ and $b_{4\beta+2,i}$ such that
\begin{align*}
L_{4\beta+2}&=5^{2\beta+1}a_{4\beta+2,0}(ut-u^2t)
+5^{2\beta+2}a_{4\beta+2,1}u+5^{2\beta+2}a_{4\beta+2,2}u^2+
5^{2\beta+1}\sum_{i=3}^{(5^{4\beta+3}-5)/3} 5^{\lfloor
\frac{5i}{6}\rfloor} a_{4\beta+2,i}u^i\\[5pt]
&\quad+5^{2\beta+2}b_{4\beta+2,1}ut+5^{2\beta+2}b_{4\beta+2,2}u^2t+
5^{2\beta+1}\sum_{i=3}^{(5^{4\beta+3}-5)/3} 5^{\lfloor
\frac{5i}{6}\rfloor} b_{4\beta+2,i}u^it,
\end{align*}
which is exactly \eqref{L-4n-2}.
Multiplying by $A$ on both sides of the above identity, and then acting
the $U_5$-operator gives
\begin{align*}
L_{4\beta+3}&=U_5(AL_{4\beta+2})\notag\\[5pt]
&=5^{2\beta+1}a_{4\beta+2,0}U_5(A(ut-u^2t))
+5^{2\beta+2}a_{4\beta+2,1}U_5(Au)
+5^{2\beta+2}a_{4\beta+2,2}U_5(Au^2)\notag\\[5pt]
&\quad+
5^{2\beta+1}\sum_{i=3}^{(5^{4\beta+3}-5)/3} 5^{\lfloor
\frac{5i}{6}\rfloor} a_{4\beta+2,i}U_5(Au^i)+5^{2\beta+2}b_{4\beta+2,1}U_5(Aut)
+5^{2\beta+2}b_{4\beta+2,2}U_5(Au^2t)\notag\\[5pt]
&\quad+
5^{2\beta+1}\sum_{i=3}^{(5^{4\beta+3}-5)/3} 5^{\lfloor
\frac{5i}{6}\rfloor} b_{4\beta+2,i}U_5(Au^it).
\end{align*}
When $i\geq 8$, $\left\lceil \frac{i+3}{5}\right\rceil\geq 3$,
and
\begin{align*}
\bigg\lfloor \frac{5i-6}{6}\bigg\rfloor>\bigg\lfloor \frac{5(\lceil \frac{i+3}{5} \rceil)+1}{6}\bigg\rfloor.
\end{align*}
For any $i\geq 3$ and $j$,
\begin{align*}
\bigg\lfloor \frac{5i-6}{6}\bigg\rfloor
+\bigg\lfloor \frac{5j-i+1}{6}\bigg\rfloor
\geq \bigg\lfloor \frac{5j+4i-11}{6}\bigg\rfloor
>  \bigg\lfloor \frac{5j}{6}\bigg\rfloor.
\end{align*}
By means of \eqref{U-ut-u2t} and Lemmas \ref{U-A-u-n} and \ref{U-A-u-n-t},
there exist integers $a_{4\beta+3,i}$ and $b_{4\beta+3,i}$ such that
\begin{align*}
L_{4\beta+3}&=5^{2\beta+1}a_{4\beta+3,0}(7u-ut+4u^2t-15u^2)
+5^{2\beta+3}a_{4\beta+3,1}u+5^{2\beta+3}a_{4\beta+3,2}u^2\notag\\[5pt]
&\quad+
5^{2\beta+2}\sum_{i=3}^{(5^{4\beta+4}-1)/3} 5^{\lfloor
\frac{5i}{6}\rfloor} a_{4\beta+3,i}u^i
+5^{2\beta+3}b_{4\beta+3,1}ut+5^{2\beta+3}b_{4\beta+3,2}u^2t
+5^{2\beta+3}b_{4\beta+3,3}u^3t\notag\\[5pt]
&\quad+
5^{2\beta+2}\sum_{i=4}^{(5^{4\beta+4}-1)/3} 5^{\lfloor
\frac{5i}{6}\rfloor} b_{4\beta+3,i}u^it,
\end{align*}
which is \eqref{L-4n-3}.
Applying the $U_5$-operator to the above identity,
we obtain
\begin{align*}
L_{4\beta+4}&=U_5(L_{4\beta+3})\\[5pt]
&=5^{2\beta+1}a_{4\beta+3,0}U_5(7u-ut+4u^2t-15u^2)
+5^{2\beta+3}a_{4\beta+3,1}U_5(u)+5^{2\beta+3}a_{4\beta+3,2}U_5(u^2)\notag\\[5pt]
&\quad+
5^{2\beta+2}\sum_{i=3}^{(5^{4\beta+4}-1)/3} 5^{\lfloor
\frac{5i}{6}\rfloor} a_{4\beta+3,i}U_5(u^i)
+5^{2\beta+3}b_{4\beta+3,1}U_5(ut)
+5^{2\beta+3}b_{4\beta+3,2}U_5(u^2t)\notag\\[5pt]
&\quad
+5^{2\beta+3}b_{4\beta+3,3}U_5(u^3t)+
5^{2\beta+2}\sum_{i=4}^{(5^{4\beta+4}-1)/3} 5^{\lfloor
\frac{5i}{6}\rfloor} b_{4\beta+3,i}U_5(u^it).
\end{align*}
Since for any $i\geq 3$,
\begin{align*}
\bigg\lfloor \frac{5i}{6}\bigg\rfloor
\geq \bigg\lfloor \frac{5\lceil \frac{i+3}{5}\rceil}{6}\bigg\rfloor.
\end{align*}
Hence, it follows from \eqref{U-7-u} and Lemmas \ref{U-u-n} and
\ref{U-u-n-t} that \eqref{L-4n-4} holds.
Hence,
\begin{align}
L_{4\beta+5}&=U_5(AL_{4\beta+4})\notag\\[5pt]
&=5^{2\beta+2}a_{4\beta+4,0}U_5(Au-Au^2t)
+5^{2\beta+3}a_{4\beta+4,1}U_5(Au)+5^{2\beta+3}a_{4\beta+4,2}U_5(Au^2)\notag\\[5pt]
&\quad+
5^{2\beta+2}\sum_{i=3}^{(5^{4\beta+5}-5)/3} 5^{\lfloor
\frac{5i}{6}\rfloor} a_{4\beta+4,i}U_5(Au^i)+5^{2\beta+3}b_{4\beta+4,1}U_5(Aut)
+5^{2\beta+3}b_{4\beta+4,2}U_5(Au^2t)\notag\\[5pt]
&\quad+
5^{2\beta+2}\sum_{i=3}^{(5^{4\beta+5}-5)/3} 5^{\lfloor
\frac{5i}{6}\rfloor} b_{4\beta+4,i}U_5(Au^it).\label{L-4-n-5}
\end{align}
By Lemma \ref{U-A-u-n},
we have
\begin{align*}
&\quad 5^{2\beta+2}\sum_{i=3}^{(5^{4\beta+5}-5)/3} 5^{\lfloor
\frac{5i}{6}\rfloor} a_{4\beta+4,i}U_5(Au^i)\notag\\[5pt]
&=5^{2\beta+4}a_{4\beta+4,3}\bigg(c_{3,1} u+c_{3,2} u^{2}
+\sum_{i=3}^{23}
5^{\lfloor \frac{5i-2}{6}\rfloor} c_{3,i}u^i
+d_{3,2} u^{2} t +5 d_{3,3} u^{3} t+\sum_{i=4}^{23}
5^{\lfloor \frac{5i-2}{6}\rfloor} d_{3,i}u^it\bigg)\notag\\[5pt]
&\quad+5^{2\beta+5}a_{4\beta+4,4}\bigg(c_{4,1} u+c_{4,2} u^{2}
+\sum_{i=3}^{28}
5^{\lfloor \frac{5i-3}{6}\rfloor} c_{4,i}u^i
+d_{4,2} u^{2} t +5 d_{4,3} u^{3} t+\sum_{i=4}^{28}
5^{\lfloor \frac{5i-3}{6}\rfloor} d_{4,i}u^it\bigg)\notag\\[5pt]
&\quad+5^{2\beta+2}\sum_{i=5}^{(5^{4\beta+5}-5)/3}
5^{\lfloor
\frac{5i}{6}\rfloor} a_{4\beta+4,i}
\bigg(c_{i,\lceil \frac{i+1}{5}\rceil} u^{\lceil \frac{i+1}{5}\rceil}+c_{i,\lceil \frac{i+1}{5}\rceil+1} u^{\lceil \frac{i+1}{5}\rceil+1}
+\sum_{j=\lceil \frac{i+1}{5}\rceil+2}^{5i+8}
5^{\lfloor \frac{5j-i+1}{6}\rfloor} c_{i,j}u^j
\notag\\[5pt]
&\quad+d_{i,\lceil \frac{i+3}{5}\rceil} u^{\lceil \frac{i+3}{5}\rceil} t + d_{i,\lceil \frac{i+3}{5}\rceil+1} u^{\lceil \frac{i+3}{5}\rceil+1} t+\sum_{j=\lceil \frac{i+3}{5}\rceil+2}^{5i+8}
5^{\lfloor \frac{5j-i+1}{6}\rfloor} d_{i,j}u^jt\bigg),
\end{align*}
where $c_{i,j}$ and $d_{i,j}$ are integers.
Since for $i\geq 5$, $\lceil \frac{i+3}{5}\rceil+1\geq 3$,
\begin{align*}
\left\lfloor \frac{5i-6}{6}\right\rfloor>
\left\lfloor \frac{5(\big\lceil \frac{i+3}{5}\big\rceil+1)}{6}\right\rfloor
\geq \left\lfloor \frac{5(\big\lceil \frac{i+1}{5}\big\rceil+1)}{6}\right\rfloor
\geq \left\lfloor \frac{5\big\lceil \frac{i+3}{5}\big\rceil}{6}\right\rfloor
\geq \left\lfloor \frac{5\big\lceil \frac{i+1}{5}\big\rceil}{6}\right\rfloor,
\end{align*}
and for any $j$,
\begin{align*}
\left\lfloor \frac{5i-6}{6}\right\rfloor+
\left\lfloor \frac{5j-i+1}{6}\right\rfloor
\geq \left\lfloor \frac{5j+4i-11}{6}\right\rfloor
>\left\lfloor \frac{5j}{6}\right\rfloor,
\end{align*}
there exist integers $c_j$ and $d_j$ such that
\begin{align}
&  5^{2\beta+2}\sum_{i=3}^{(5^{4\beta+5}-5)/3} 5^{\lfloor
\frac{5i}{6}\rfloor} a_{4\beta+4,i}U_5(Au^i)
=5^{2\beta+4}c_{1}
u+5^{2\beta+4}c_{2} u^{2}
 \notag\\[5pt]
&\qquad \qquad +5^{2\beta+3}\sum_{j=3}^{(5^{4\beta+6}-1)/3}
5^{\lfloor \frac{5j}{6}\rfloor} c_ju^j +5^{2\beta+4}d_{2} u^{2} t
+5^{2\beta+3}\sum_{j=3}^{(5^{4\beta+6}-1)/3} 5^{\lfloor
\frac{5j}{6}\rfloor} d_ju^jt.  \label{U-A-u-comp}
\end{align}
Analogously, by Lemma \ref{U-A-u-n-t}, we obtain that
there exist integers $c_j'$ and $d_j'$ such that
\begin{align}
&  5^{2\beta+2}\sum_{i=3}^{(5^{4\beta+5}-5)/3} 5^{\lfloor
\frac{5i}{6}\rfloor} b_{4\beta+4,i}U_5(Au^it) =5^{2\beta+4}c_{1}'
u+5^{2\beta+4}c_{2}' u^{2}
\notag  \\[5pt]
&\qquad \qquad+5^{2\beta+3}\sum_{j=3}^{(5^{4\beta+6}-16)/3}
5^{\lfloor \frac{5j}{6}\rfloor} c_j'u^j +5^{2\beta+4}d_{2}' u^{2} t
+5^{2\beta+3}\sum_{j=3}^{(5^{4\beta+6}-16)/3} 5^{\lfloor
\frac{5j}{6}\rfloor} d_j'u^jt. \label{U-A-u-t-comp}
\end{align}
Substituting \eqref{U-A-u-comp}, \eqref{U-A-u-t-comp},
\eqref{A-u-f}, \eqref{A-u-2-f},
\eqref{A-u-t-f}, \eqref{A-u-2-t-f} and \eqref{U-u-u2}
into \eqref{L-4-n-5} gives that \eqref{L-4n-1}
is also true for $\beta+1$.

Therefore, by induction, \eqref{L-4n-1}--\eqref{L-4n-4}
hold for all $\beta\geq 0$.
\qed

\subsection{Proof of Theorem \ref{Th-1}}\label{Sect-pfthm}

In this subsection, we give a proof of Theorem \ref{Th-1}.

\noindent{\it Proof  of Theorem \ref{Th-1}}.
The congruence \eqref{1-7} is an immediate result of
\eqref{L-even-gf}, \eqref{L-4n-2} and \eqref{L-4n-4}.
In order to prove \eqref{1-8},
we only need to prove that for any $\beta\geq 0$,
\begin{align}
PDN1\left(5^{4\beta+3}n +\frac{r\cdot 5^{4\beta+2}+1}{ 24}
 \right)\equiv 0 \pmod {5^{ 2\beta
 +2}},\label{cong-4-n-3}
\end{align}
and
\begin{align}
PDN1\left(5^{4\beta+5}n +\frac{r\cdot 5^{4\beta+4}+1}{ 24}
 \right)\equiv 0 \pmod {5^{2\beta
 +3}},\label{cong-4-n-5}
\end{align}
where $r\in\{71,119\}$.

We first prove \eqref{cong-4-n-3}. Thanks to
  \eqref{L-even-gf} and
\eqref{L-4n-2},
\begin{align*}
\sum_{n=0}^\infty PDN1
\left(
5^{4\beta+2}n+\frac{23\cdot 5^{4\beta+2}+1}{24}\right)q^n
\equiv
5^{2\beta+1}a_{4\beta+2,0}\frac{J_{2}^2}{qJ_1^5}(ut-u^2t) \pmod{5^{2\beta+2}}.
\end{align*}
Define
\begin{align*}
\sum_{n=0}^\infty E(n)q^n&:=\frac{J_{2}^2}{qJ_1^5}(ut-u^2t)
=\frac{J_{10}^4}{J_1^3J_2^2J_5^2}-q^2
\frac{J_{10}^{11}}{J_1^4J_2^5J_5^5}.
\end{align*}
Thus, \eqref{cong-4-n-3} is equivalent to
\begin{align}
E(5n+t)\equiv 0\pmod{5},\label{A-5-cong}
\end{align}
where $t=2,4$.
By Theorem \ref{con_F_modular_function}, we find
\begin{align*}
F_1(\tau)=
q^{\frac{71}{120}}\frac{\eta(5\tau)\eta^2(10\tau)\eta_{10,4}^4(\tau)}{\eta^2_{10,5}(\tau)}
\sum_{n=0}^\infty E(5n+2)q^n
\end{align*}
and
\begin{align*}
F_2(\tau)=
q^{\frac{119}{120}}\frac{\eta(5\tau)\eta^2(10\tau)\eta_{10,5}^6(\tau)}{\eta^4_{10,4}(\tau)}
\sum_{n=0}^\infty E(5n+4)q^n
\end{align*}
are modular functions for $\Gamma_1(10)$.
The set of inequivalent cusps of $\Gamma_1(10)$ is given by
\begin{align*}
\left\{0,\frac{1}{5},\frac{1}{4},\frac{3}{10},
\frac{1}{3},\frac{2}{5},\frac{1}{2},\infty\right\}.
\end{align*}
Let
\begin{align*}
F_1'(\tau)=\eta^{32}(\tau)\eta^8(2\tau)F_1(\tau)
\end{align*}
and
\begin{align*}
F_2'(\tau)=\eta^{24}(\tau)\eta^{24}(2\tau)F_2(\tau).
\end{align*}
By Theorem \ref{order_F-1},
we obtain that orders of $F_1'(\tau)$ and $F_2'(\tau)$
are non-negative at each cusps of $\Gamma_1(10)$.
Hence, $F_1'(\tau)$ and $F_2'(\tau)$
are modular forms with weight 20 and 24 for $\Gamma_1(10)$,
respectively.
By the $q$-expansions of $F_1'(\tau)$ and $F_2'(\tau)$,
one can check the coefficients of their first 150 terms
are congruent to 0 modulo 5.
It follows from Sturm's theorem \ref{Sturm-theorem} that
$F_1'(\tau)\equiv 0\pmod{5}$ and
$F_2'(\tau)\equiv 0\pmod{5}$,
which implies \eqref{A-5-cong}.

Now, we prove \eqref{cong-4-n-5}.
Combining \eqref{L-even-gf} and \eqref{L-4n-4} yields
\begin{align*}
\sum_{n=0}^\infty PDN1
\left(
5^{4\beta+4}n+\frac{23\cdot 5^{4\beta+4}+1}{24}\right)q^n
\equiv
5^{2\beta+2}a_{4\beta+4,0}\frac{J_{2}^2}{qJ_1^5}(u-u^2t) \pmod{5^{2\beta+3}}.
\end{align*}
Consider
\begin{align*}
\sum_{n=1}^\infty G(n)q^n&:=\frac{J_{2}^2}{qJ_1^5}(u-u^2t)
=q\frac{J_{10}^7}{J_1^6J_2J_5^3}-q^2
\frac{J_{10}^{11}}{J_1^4J_2^5J_5^5}.
\end{align*}
Thus, \eqref{cong-4-n-5} is equivalent to
\begin{align}
G(5n+t)\equiv 0\pmod{5},\label{B-5-cong}
\end{align}
where $t=2,4$.
By Theorem \ref{con_F_modular_function}, we find
\begin{align*}
H_1(\tau)=
q^{\frac{71}{120}}\frac{\eta(5\tau)\eta^2(10\tau)\eta_{10,4}^4(\tau)}{\eta^2_{10,5}(\tau)}
\sum_{n=0}^\infty G(5n+2)q^n
\end{align*}
and
\begin{align*}
H_2(\tau)=
q^{\frac{119}{120}}\frac{\eta(5\tau)\eta^2(10\tau)\eta_{10,5}^6(\tau)}{\eta^4_{10,4}(\tau)}
\sum_{n=0}^\infty G(5n+4)q^n
\end{align*}
are modular functions for $\Gamma_1(10)$.
Let
\begin{align*}
H_1'(\tau)=\eta^{32}(\tau)\eta^8(2\tau)H_1(\tau)
\end{align*}
and
\begin{align*}
H_2'(\tau)=\eta^{24}(\tau)\eta^{24}(2\tau)H_2(\tau).
\end{align*}
Using Theorem \ref{order_F-1} to compute the lower bounds of orders
of $H_1'(\tau)$ and $H_2'(\tau)$, we find $H_1'(\tau)$ and
$H_2'(\tau)$ are modular forms with weight 20 and 24 for
$\Gamma_1(10)$, respectively. By the $q$-expansions of $H_1'(\tau)$
and $H_2'(\tau)$, one can check the coefficients of their first 150
terms are congruent to 0 modulo 5. It follows from Sturm's theorem
\ref{Sturm-theorem} that $H_1'(\tau)\equiv 0\pmod{5}$ and
$H_2'(\tau)\equiv 0\pmod{5}$, which implies \eqref{B-5-cong}.

Therefore, \eqref{1-8} holds for any $\alpha\geq 1$.
\qed

\section{Congruences modulo the powers of $7$}\label{sect-pf-th2}

\subsection{The generating functions}

To prove Theorem \ref{cong}, we derive the expressions
of the generating functions of
infinite arithmetic progressions on $PDN1(n)$ in \eqref{1-9} and
\eqref{1-10}.

Let
\begin{align*}
A:=\frac{\eta^2(2 \tau ) \eta^5(49\tau
)}{\eta^5(\tau)\eta^2(98\tau)}.
\end{align*}
Set $L_0=1$ and for $\alpha\geq 1$, define
\begin{align*}
L_{2\alpha-1}:=U_7(AL_{2\alpha-2}), \quad
L_{2\alpha}:=U_7(L_{2\alpha-1}).
\end{align*}
It is easy to see from \eqref{1-1} that
\begin{align*}
L_{2\alpha-1}&=\frac{J_7^5}{J_{14}^2}q\sum_{n=0}^\infty
PDN1\left( 7^{2\alpha-1}n+\frac{17\cdot 7^{2\alpha-1}+1}{24}\right)q^n,\\[5pt]
L_{2\alpha}&=
\frac{J_1^5}{J_{2}^2}q\sum_{n=0}^\infty PDN1\left(
7^{2\alpha}n+\frac{23\cdot 7^{2\alpha}+1}{24}\right)q^n.
\end{align*}
It is obvious from Theorem \ref{GHN-THM} that
$A$ is a modular function for $\Gamma_0(98)$.
Thus, $U_7(A)$ is a modular function for $\Gamma_0(14)$,
and, by Theorem \ref{THM-cusp} and Lemma \ref{GH-ord}, we have
\begin{align}
\mathrm{ord}_0(U_7(A)) \geq -16,\quad \mathrm{ord}_{1/2}(U_7(A))\geq -2, \quad \mathrm{ord}_{1/7}(U_7(A)) \geq 3,
\quad \mathrm{ord}_\infty(U_7(A)) \geq 1.\label{order-U-7-A}
\end{align}
Hence, $L_\alpha (\alpha\geq 1)$ are modular functions
for $\Gamma_0(14)$.
Define
\begin{align*}
t:=\frac{\eta(2\tau)
\eta^7(7\tau)}{\eta(\tau)\eta^7(14\tau)},
\quad z:=\frac{\eta^3(\tau)\eta(2\tau)\eta^3(7\tau)}
{\eta^7(14 \tau)},
 \quad u:=z^{-1}=\frac{\eta^7(14
\tau)}{\eta^3(\tau)\eta(2\tau)\eta^3(7\tau)}.
\end{align*}
The functions $t, z$ and $u$ are modular functions for $\Gamma_0(14)$,
and the orders of $t$ and $z$ at the cusps of $\Gamma_0(14)$ are
\begin{align}
\mathrm{ord}_0(t) &=0, &\mathrm{ord}_{1/2}(t) &=0, &\mathrm{ord}_{1/7}(t) &=2, &\mathrm{ord}_\infty(t) &=-2,\label{order-t-7}\\
\mathrm{ord}_0(z) &=2, &\mathrm{ord}_{1/2}(z) &=1, &\mathrm{ord}_{1/7}(z) &=0, &\mathrm{ord}_\infty(z) &=-3.\label{order-z-7}
\end{align}
Hence, $z^8U_7(A)$ has a pole only at $\infty$.
By Algorithm AB of Radu \cite{Radu},
\begin{align*}
\langle E^\infty(14) \rangle_\mathbb{Q}
=\mathbb{Q}[t]+z\mathbb{Q}[t].
\end{align*}
Using the algorithm of Frye and Garvan \cite{FG} or
the zero recognition method of Paule and
Radu \cite{PR21}, we have
\begin{align}
t^3&=(z^{2}-16z+64)+(9z-80)t +17 t^2,\label{t-3}\\[5pt]
t^4 &=(17z^{2}-272z+1088)+ (z^{2}+137z-1296)t+
(9z+209)t^2.\label{t-4}
\end{align}
Therefore, $\langle E^\infty(14) \rangle_\mathbb{Q}$
can be expressed as
\begin{align*}
\langle E^\infty(14) \rangle_\mathbb{Q}
=\mathbb{Q}[z]+t\mathbb{Q}[z]+t^2\mathbb{Q}[z].
\end{align*}

Theorem \ref{cong} is an immediate result of the following theorem,
which gives the expression of $L_{\alpha}$ for $\alpha\geq 1$.

\begin{theorem}\label{L-expression-7}
For any $\alpha\geq 1$, there exist integers $a_{\alpha,i},
b_{\alpha,i},c_{\alpha,i}$ such that
\begin{align}
L_{2\alpha-1}&=7^{\alpha+1}a_{2\alpha-1,1}u
+7^{\alpha+1}a_{2\alpha-1,2}u^2+
7^{\alpha+1}\sum_{i=3}^{(7^{2\alpha}-1)/6} 7^{\lceil
\frac{7i-24}{26}\rceil} a_{2\alpha-1,i}u^i
 +7^{\alpha}b_{2\alpha-1,1}ut\notag\\[5pt]
&\quad
 +7^{\alpha}b_{2\alpha-1,2}u^2t+
7^{\alpha+1}\sum_{i=3}^{(7^{2\alpha}-1)/6} 7^{\lceil
\frac{7i-24}{26}\rceil} b_{2\alpha-1,i}u^it +
7^{\alpha}c_{2\alpha-1,2}u^2t^2\notag\\[5pt]
&\quad+ 7^{\alpha+1}\sum_{i=3}^{(7^{2\alpha}-1)/6} 7^{\lceil
\frac{7i-24}{26}\rceil}c_{2\alpha-1,i}
u^it^2,\label{L-2n-1}\\[5pt]
L_{2\alpha}&=7^{\alpha+1}a_{2\alpha,1}u
+7^{\alpha+1}a_{2\alpha,2}u^2+
7^{\alpha+1}\sum_{i=3}^{(7^{2\alpha+1}-7)/6} 7^{\lceil
\frac{7i+3}{26}\rceil}a_{2\alpha,i}u^i
 +7^{\alpha+1}b_{2\alpha,1}ut\notag\\[5pt]
&\quad
 +7^{\alpha+1}b_{2\alpha,2}u^2t+
7^{\alpha+1}\sum_{i=3}^{(7^{2\alpha+1}-7)/6} 7^{\lceil
\frac{7i+3}{26}\rceil}b_{2\alpha,i}u^it+
7^{\alpha+1}c_{2\alpha,2}u^2t^2\notag\\[5pt]
&\quad + 7^{\alpha+1}\sum_{i=3}^{(7^{2\alpha+1}-7)/6} 7^{\lceil
\frac{7i+3}{26}\rceil}c_{2\alpha,i} u^it^2.\label{L-2n}
\end{align}
\end{theorem}

\subsection{Fundamental lemmas}

To prove Theorem \ref{L-expression-7}, we require the following
modular equation  of the seventh order.

\begin{lemma}\label{modular-equation-7}
Let $\sigma_i(q)~(1\leq i\leq 7)$ be
 defined in the Appendix B. Then
\begin{align}\label{mod-eq-7}
u^7=\sigma_1(q^7) u^6+\sigma_2(q^7) u^5+\sigma_3(q^7) u^4
+\sigma_4(q^7) u^3 +\sigma_5(q^7) u^2+\sigma_6(q^7) u+\sigma
_7(q^7).
\end{align}
\end{lemma}

\pf
Since $u, t$ are modular functions for $\Gamma_0(14)$,
as well as for
$\Gamma_0(98)$, \eqref{mod-eq-7} can be shown by utilizing the algorithm of Frye and Garvan \cite{FG} or
the zero recognition method of Paule and
Radu \cite{PR21}.
\qed

\begin{corollary}
For any function $f\colon \mathbb{H}\rightarrow \mathbb{C}$
and any $n\geq 7$, we have
\begin{align}
U_7(fu^n)=\sum_{i=1}^7 \sigma_i(q) U_7(fu^{n-i}).\label{u-f-u-7}
\end{align}
\end{corollary}

Utilizing the modular equation \eqref{mod-eq-7}, we deduce
  the following lemma.

\begin{lemma}
Let $v_0(n)=\left\lceil \frac{n}{7}\right\rceil$,
$v_1(n)=\left\lceil \frac{n+5}{7}\right\rceil$, $v_2(n)=\left\lceil
\frac{n+3}{7}\right\rceil+1$,
$\mu(n,i)=\left\lceil\frac{7i-n+6}{26}\right\rceil$. Then for any
$n\geq 1$, there exist integers $d_{n,i}, e_{n,i}$ and $h_{n,i}$
such that
\begin{align}
U_7(u^n)&=
\begin{cases}
& d_{n,v_0(n)}u^{v_0(n)} +\sum\limits_{i=v_0(n)+1}^{7n} 7^{\mu(n,i)}
d_{n,i}u^i+ e_{n,v_1(n)}u^{v_1(n)}t
\\[5pt]
&\quad  +\sum\limits_{i=v_1(n)+1}^{7n} 7^{\mu(n,i)} e_{n,i}u^it+
\sum\limits_{i=v_2(n)}^{7n} 7^{\mu(n,i)}
h_{n,i}u^it^2, \qquad \text{ if } n\equiv 1,2 \pmod{7},\\[5pt]
& d_{n,v_0(n)}u^{v_0(n)} +\sum\limits_{i=v_0(n)+1}^{7n} 7^{\mu(n,i)}
d_{n,i}u^i
\\[5pt]
&\quad  +\sum\limits_{i=v_1(n)}^{7n} 7^{\mu(n,i)} e_{n,i}u^it+
\sum\limits_{i=v_2(n)}^{7n} 7^{\mu(n,i)} h_{n,i}u^it^2, \qquad \
 \text{ if } n\equiv 0,3,5,6 \pmod{7},\\[5pt]
& d_{n,v_0(n)}u^{v_0(n)} +d_{n,v_0(n)+1}u^{v_0(n)+1}
+\sum\limits_{i=v_0(n)+2}^{7n} 7^{\mu(n,i)} d_{n,i}u^i+
e_{n,v_1(n)}u^{v_1(n)}t
\\[5pt]
&\quad  +\sum\limits_{i=v_1(n)+1}^{7n} 7^{\mu(n,i)} e_{n,i}u^it+
h_{n,v_2(n)}u^{v_2(n)}t^2\\[5pt]
&\quad+ \sum\limits_{i=v_2(n)+1}^{7n}
 7^{\mu(n,i)} h_{n,i}u^it^2, \qquad \qquad
 \qquad \qquad \qquad \qquad
\text{ if } n\equiv 4 \pmod{7}.
\end{cases}\label{U-7-u-n}
\end{align}
\end{lemma}

\pf We prove \eqref{U-7-u-n}
by induction on $n$. For $1\leq n \leq 7$,
we get \eqref{U-7-u-n} holds from our Mathematica supplement.
Suppose that \eqref{U-7-u-n} is true for all
 integers less than $n$ with $n\geq
8$. Now we prove that it holds for $n$ according to the residues of
$n$ modulo 7. We only prove the case $n\equiv 0\pmod{7}$, the rest
can be shown analogously. Let $n=7 k$, where $k\geq 2$. Then,
by \eqref{u-f-u-7},
\begin{align*}
U_7(u^n)&=\sum_{i=1}^7 \sigma_i(q) U_7(u^{n-i}).
\end{align*}

We prove each term can be expressed in the form of
\begin{align}
a_{k}'u^k+\sum_{i=k+1}^{49k}7^{\lceil\frac{7i-7k+6}{26}\rceil}a'_i
u^i +\sum_{i=k+1}^{49k}7^{\lceil\frac{7i-7k+6}{26}\rceil}b'_i
u^it +\sum_{i=k+2}^{49k}7^{\lceil\frac{7i-7k+6}{26}\rceil}
c'_iu^it^2,\label{u-7k-form}
\end{align}
where $a'_i, b'_i,c'_i$ are integers. Here, we only prove that it
holds for $\sigma_7(q)U_7(u^{n-7})$,
and the other terms can be shown in the
same manner.

From the Appendix B, we have
\[\sigma_7(q)=A_0+A_1t+A_2t^2,\]
where
\begin{align*}
A_0&:= u + 6818 u^2 + 1384397 u^3 + 17815420 u^4 - 419183387 u^5 +
 311299254 u^6 + 9321683217 u^7,
 \\[6pt]
 A_1& :=42 u^2 + 96285 u^3 + 14605283 u^4 + 232239126 u^5 - 2075328360 u^6 -
 10734059462 u^7,\\[6pt]
 A_2&:=
 735 u^3 + 521017 u^4 + 42824236 u^5 + 634128110 u^6 + 1412376245
 u^7.
\end{align*}
 Then, by the inductive hypothesis,
\begin{align*}
&\quad \sigma_7(q)U_7(u^{n-7})\\[5pt]
&=\sigma_7(q)\bigg(d_{n-7,k-1}u^{k-1}+\sum_{i=k}^{49k-49}
7^{\lceil\frac{7i-7k+13}{26}\rceil} d_{n-7,i}u^i
\nonumber\\[6pt]
&\quad+\sum_{i=k}^{49k-49} 7^{\lceil\frac{7i-7k+13}{26}\rceil}
e_{n-7,i}u^it +\sum_{i=k+1}^{49k-49}
7^{\lceil\frac{7i-7k+13}{26}\rceil}
 h_{n-7,i}u^it^2\bigg)\\[5pt]
&=A_0\bigg(d_{n-7,k-1}u^{k-1}+\sum_{i=k}^{49k-49}
7^{\lceil\frac{7i-7k+13}{26}\rceil}
 d_{n-7,i}u^i\bigg) \\[5pt]
&\quad+\bigg(A_1 \bigg(d_{n-7,k-1}u^{k-1} +\sum_{i=k}^{49k-49}
7^{\lceil\frac{7i-7k+13}{26}\rceil}
 d_{n-7,i}u^i\bigg)
  +A_0 \sum_{i=k}^{49k-49} 7^{\lceil\frac{7i-7k+13}{26}\rceil}
e_{n-7,i}u^i\bigg)t   \nonumber\\[5pt]
&\quad +\bigg(A_0 \sum_{i=k+1}^{49k-49}
7^{\lceil\frac{7i-7k+13}{26}\rceil} h_{n-7,i}u^i +A_1
\sum_{i=k}^{49k-49} 7^{\lceil\frac{7i-7k+13}{26}\rceil}
 e_{n-7,i}u^i  \nonumber\\[5pt]
&\quad +A_2 \bigg(d_{n-7,k-1}u^{k-1} +\sum_{i=k}^{49k-49}
7^{\lceil\frac{7i-7k+13}{26}\rceil}
 d_{n-7,i}u^i\bigg)\bigg)t^2  +\bigg(A_2\sum_{i=k}^{49k-49}
7^{\lceil\frac{7i-7k+13}{26}\rceil} e_{n-7,i}u^i
\\[5pt]
&\quad+A_1\sum_{i=k+1}^{49k-49} 7^{\lceil\frac{7i-7k+13}{26}\rceil}
 h_{n-7,i}u^i\bigg)t^3 +A_2 \sum_{i=k+1}^{49k-49}
7^{\lceil\frac{7i-7k+13}{26}\rceil}
 h_{n-7,i}u^it^4\\[5pt]
 &=
d_{n-7,k-1}u^k+\sum_{i=k+1}^{49k-42}7^{\lceil\frac{7i-7k+6}{26}\rceil}
d_{i}u^i +\sum_{i=k+1}^{49k-42}7^{\lceil\frac{7i-7k+6}{26}\rceil}
e_{i}u^it +\sum_{i=k+2}^{49k-42}7^{\lceil\frac{7i-7k+6}{26}\rceil}
f_{i}u^it^2\\[5pt]
&\quad +\sum_{i=k+3}^{49k-42}7^{\lceil\frac{7i-7k+6}{26}\rceil}
g_{i}u^it^3 +\sum_{i=k+4}^{49k-42}7^{\lceil\frac{7i-7k+6}{26}\rceil}
h_{i}u^it^4,
\end{align*}where $d_i, e_i, f_i, g_i, h_i$ are integers.
After substituting \eqref{t-3} and \eqref{t-4} into the above equation,
we obtain that there exist
integers $a_i', b_i', c_i'$ such that
\begin{align*}
\sigma_7(q)U_7(u^{n-7})=
d_{n-7,k-1}u^k+\sum_{i=k+1}^{49k-42}7^{\lceil\frac{7i-7k+6}{26}\rceil}
a_i'u^i +\sum_{i=k+1}^{49k-42}7^{\lceil\frac{7i-7k+6}{26}\rceil}
b_i'u^it +\sum_{i=k+2}^{49k-42}7^{\lceil\frac{7i-7k+6}{26}\rceil}
c_i'u^it^2,
\end{align*}
that is, $\sigma_7(q)U_7(u^{n-7})$ has the form of \eqref{u-7k-form}.

Therefore, \eqref{U-7-u-n} is true for $n\equiv 0\pmod{7}$.
\qed

Analogously, we obtain the expressions of $U_7(u^nt), U_7(u^nt^2),
U_7(Au^n), U_7(Au^nt), U_7(Au^nt^2)$.

\begin{lemma}\label{U-7-ut}
Let $v_0(n)=\left\lceil \frac{n}{7}\right\rceil$,
$v_1(n)=\left\lceil \frac{n+4}{7}\right\rceil$, $v_2(n)=\left\lceil
\frac{n+2}{7}\right\rceil+1$,
$\mu(n,i)=\left\lceil\frac{7i-n+6}{26}\right\rceil$. Then for $n\geq
1$, there exist integers $d_{n,i}, e_{n,i}$ and $h_{n,i}$ such that
%
%
%
\begin{align*}
U_7(u^nt)&=
\begin{cases}
& \sum\limits_{i=v_0(n)}^{7n} 7^{\mu(n,i)}
d_{n,i}u^i+\sum\limits_{i=v_1(n)}^{7n} 7^{\mu(n,i)}
e_{n,i}u^it\\[5pt]
&\quad+ \sum\limits_{i=v_2(n)}^{7n} 7^{\mu(n,i)}
h_{n,i}u^it^2,\qquad \qquad\qquad \qquad\qquad\qquad \qquad \text{ if } n\equiv 1,2,6 \pmod{7},\\[5pt]
& 7d_{n,v_0(n)}u^{v_0(n)} +\sum\limits_{i=v_0(n)+1}^{7n}
7^{\mu(n,i)} d_{n,i}u^i+ e_{n,v_1(n)}u^{v_1(n)}t
\\[5pt]
&\quad +\sum\limits_{i=v_1(n)+1}^{7n} 7^{\mu(n,i)} e_{n,i}u^it+
\sum\limits_{i=v_2(n)}^{7n} 7^{\mu(n,i)}
h_{n,i}u^it^2,\qquad   \qquad  \text{ if } n\equiv 3 \pmod{7},\\[5pt]
& d_{n,v_0(n)}u^{v_0(n)} +d_{n,v_0(n)+1}u^{v_0(n)+1}
+\sum\limits_{i=v_0(n)+2}^{7n} 7^{\mu(n,i)}
d_{n,i}u^i+e_{n,v_1(n)}u^{v_1(n)}t
\\[5pt]
&\quad +\sum\limits_{i=v_1(n)+1}^{7n} 7^{\mu(n,i)} e_{n,i}u^it
+h_{n,v_2(n)}u^{v_2(n)}t^2\\[5pt]
&\quad + \sum\limits_{i=v_2(n)+1}^{7n} 7^{\mu(n,i)}
h_{n,i}u^it^2, \qquad \qquad\qquad \quad \
\qquad\qquad \qquad \text{ if } n\equiv 4,5 \pmod{7},\\[5pt]
& d_{n,v_0(n)}u^{v_0(n)} +d_{n,v_0(n)+1}u^{v_0(n)+1}
+\sum\limits_{i=v_0(n)+2}^{7n} 7^{\mu(n,i)} d_{n,i}u^i
 +e_{n,v_1(n)}u^{v_1(n)}t\\[5pt]
&\quad +\sum\limits_{i=v_1(n)+1}^{7n} 7^{\mu(n,i)} e_{n,i}u^it +
\sum\limits_{i=v_2(n)}^{7n} 7^{\mu(n,i)}
h_{n,i}u^it^2, \qquad   \qquad\text{ if }
n\equiv 0 \pmod{7}.
\end{cases}
\end{align*}
\end{lemma}

\begin{lemma}
Let $v_0(n)=\left\lceil \frac{n-1}{7}\right\rceil$,
$v_1(n)=\left\lceil \frac{n+4}{7}\right\rceil$, $v_2(n)=\left\lceil
\frac{n+1}{7}\right\rceil+1$,
$\mu(n,i)=\left\lceil\frac{7i-n+6}{26}\right\rceil$. Then for $n\geq
1$, there exist integers $d_{n,i}, e_{n,i}$ and $h_{n,i}$ such that
\begin{align*}
U_7(u^nt^2)&=
\begin{cases}
& d_{n,v_0(n)}u^{v_0(n)} +d_{n,v_0(n)+1}u^{v_0(n)+1}
+\sum\limits_{i=v_0(n)+2}^{7n} 7^{\mu(n,i)} d_{n,i}u^i+
e_{n,v_1(n)}u^{v_1(n)}t
\\[5pt]
&\quad +\sum\limits_{i=v_1(n)+1}^{7n} 7^{\mu(n,i)} e_{n,i}u^it+
\sum\limits_{i=v_2(n)}^{7n} 7^{\mu(n,i)}
h_{n,i}u^it^2,\qquad \qquad \qquad \qquad \text{ if } n\equiv 0,1 \pmod{7},\\[5pt]
& \sum\limits_{i=v_0(n)}^{7n} 7^{\mu(n,i)} d_{n,i}u^i
+\sum\limits_{i=v_1(n)}^{7n} 7^{\mu(n,i)} e_{n,i}u^it
+\sum\limits_{i=v_2(n)}^{7n} 7^{\mu(n,i)}
h_{n,i}u^it^2, \qquad  \text{ if } n\equiv 2 \pmod{7},\\[5pt]
& 7d_{n,v_0(n)}u^{v_0(n)}+ \sum\limits_{i=v_0(n)+1}^{7n}
7^{\mu(n,i)} d_{n,i}u^i+e_{n,v_1(n)}u^{v_1(n)}t
+\sum\limits_{i=v_1(n)+1}^{7n} 7^{\mu(n,i)}
e_{n,i}u^it \\[5pt]
&\quad +\sum\limits_{i=v_2(n)}^{7n} 7^{\mu(n,i)}
h_{n,i}u^it^2, \qquad\qquad\ \ \qquad\qquad \qquad\qquad \qquad\qquad \qquad \text{ if } n\equiv 3 \pmod{7},\\[5pt]
& d_{n,v_0(n)}u^{v_0(n)} +d_{n,v_0(n)+1}u^{v_0(n)+1}
+\sum\limits_{i=v_0(n)+2}^{7n} 7^{\mu(n,i)} d_{n,i}u^i+
e_{n,v_1(n)}u^{v_1(n)}t\\[5pt]
&\  +\sum\limits_{i=v_1(n)+1}^{7n} 7^{\mu(n,i)}
e_{n,i}u^it+h_{n,v_2(n)}u^{v_2(n)}t^2+ \sum\limits_{i=v_2(n)+1}^{7n} 7^{\mu(n,i)} h_{n,i}u^it^2,
 \text{ if }
n\equiv 4,5,6 \pmod{7}.
\end{cases}
\end{align*}
\end{lemma}

\begin{lemma}\label{U-7-A-u}
Let $v_0(n)=\left\lceil \frac{n+1}{7}\right\rceil$,
$v_1(n)=\left\lceil \frac{n+6}{7}\right\rceil$, $v_2(n)=\left\lceil
\frac{n+3}{7}\right\rceil+1$,
$\mu(n,i)=\left\lceil\frac{7i-n}{26}\right\rceil$. Then for $n\geq
0$, there exist integers $d_{n,i}, e_{n,i}$ and $h_{n,i}$ such that
\begin{align*}
U_7(Au^n)=
\begin{cases}
&\sum\limits_{i=v_0(n)}^{7n+8} 7^{\mu(n,i)} d_{n,i}u^i
+\sum\limits_{i=v_1(n)}^{7n+8} 7^{\mu(n,i)}
e_{n,i}u^it\\[5pt]
&\quad+ \sum\limits_{i=v_2(n)}^{7n+8} 7^{\mu(n,i)}
h_{n,i}u^it^2,\qquad\qquad\qquad\qquad \qquad \qquad \text{ if } n\equiv 0,2,3 \pmod{7},\\[5pt]
& 7d_{n,v_0(n)}u^{v_0(n)} +7d_{n,v_0(n)+1}u^{v_0(n)+1}
+\sum\limits_{i=v_0(n)+2}^{7n+8} 7^{\mu(n,i)}
d_{n,i}u^i+\\[5pt]
&\quad e_{n,v_1(n)}u^{v_1(n)}t+ e_{n,v_1(n)+1}u^{v_1(n)+1}t
+\sum\limits_{i=v_1(n)+2}^{7n+8} 7^{\mu(n,i)}
e_{n,i}u^it\\[5pt]
&\quad +h_{n,v_2(n)}u^{v_2(n)}t^2+ \sum\limits_{i=v_2(n)+1}^{7n+8}
7^{\mu(n,i)}
h_{n,i}u^it^2, \qquad\qquad \text{ if } n\equiv 1 \pmod{7},\\[5pt]
& \sum\limits_{i=v_0(n)}^{7n+8} 7^{\mu(n,i)} d_{n,i}u^i
+\sum\limits_{i=v_1(n)}^{7n+8} 7^{\mu(n,i)}
e_{n,i}u^it+h_{n,v_2(n)}u^{v_2(n)}t^2\\[5pt]
&\quad+ \sum\limits_{i=v_2(n)+1}^{7n+8} 7^{\mu(n,i)}
h_{n,i}u^it^2,\qquad\qquad\qquad\qquad\qquad
\qquad \text{ if } n\equiv 4 \pmod{7},\\[5pt]
& d_{n,v_0(n)}u^{v_0(n)}
+d_{n,v_0(n)+1}u^{v_0(n)+1}+d_{n,v_0(n)+2}u^{v_0(n)+2}
+\sum\limits_{i=v_0(n)+3}^{7n+8} 7^{\mu(n,i)}
d_{n,i}u^i+\\[5pt]
&\quad e_{n,v_1(n)}u^{v_1(n)}t +e_{n,v_1(n)+1}u^{v_1(n)+1}t
+\sum\limits_{i=v_1(n)+2}^{7n+8} 7^{\mu(n,i)}
e_{n,i}u^it\\[5pt]
&\quad +h_{n,v_2(n)}u^{v_2(n)}t^2 + \sum\limits_{i=v_2(n)+1}^{7n+8}
7^{\mu(n,i)} h_{n,i}u^it^2,  \qquad\ \qquad\text{ if } n\equiv 5,6
\pmod{7}.
\end{cases}
\end{align*}
\end{lemma}

\begin{lemma}
Let $v_0(n)=\left\lceil \frac{n}{7}
\right\rceil$,
$v_1(n)=\left\lceil \frac{n+5}{7}
\right\rceil$, $v_2(n)=\left\lceil
\frac{n+3}{7}\right\rceil+1$,
$\mu(n,i)=\left\lceil\frac{7i-n}{26}\right\rceil$. Then for
$n\geq 0$, there exist integers $d_{n,i}, e_{n,i}$ and $h_{n,i}$
such that
\begin{align*}
U_7(Au^nt)=
\begin{cases}
& d_{n,v_0(n)}u^{v_0(n)}
+d_{n,v_0(n)+1}u^{v_0(n)+1}+d_{n,v_0(n)+2}u^{v_0(n)+2}
+\sum\limits_{i=v_0(n)+3}^{7n+8} 7^{\mu(n,i)}
d_{n,i}u^i+\\[5pt]
&\quad e_{n,v_1(n)}u^{v_1(n)}t+e_{n,v_1(n)+1}u^{v_1(n)+1}t+
\sum\limits_{i=v_1(n)+2}^{7n+8} 7^{\mu(n,i)}
e_{n,i}u^it\\[5pt]
&\quad +h_{n,v_2(n)}u^{v_2(n)}t^2+ \sum\limits_{i=v_2(n)+1}^{7n+8}
7^{\mu(n,i)}
h_{n,i}u^it^2,
\qquad\qquad \quad \text{ if } n\equiv 0,5,6 \pmod{7},\\[5pt]
& 7d_{n,v_0(n)}u^{v_0(n)} +7d_{n,v_0(n)+1}u^{v_0(n)+1}
+\sum\limits_{i=v_0(n)+2}^{7n+8} 7^{\mu(n,i)}
d_{n,i}u^i\\[5pt]
&\quad+ e_{n,v_1(n)}u^{v_1(n)}t +e_{n,v_1(n)+1}u^{v_1(n)+1}t
+\sum\limits_{i=v_1(n)+2}^{7n+8} 7^{\mu(n,i)}
e_{n,i}u^it\\[5pt]
&\quad +h_{n,v_2(n)}u^{v_2(n)}t^2+ \sum\limits_{i=v_2(n)+1}^{7n+8}
7^{\mu(n,i)} h_{n,i}u^it^2, \quad
\qquad\qquad\text{ if } n\equiv 1,2
 \pmod{7},\\[5pt]
& \sum\limits_{i=v_0(n)}^{7n+8} 7^{\mu(n,i)} d_{n,i}u^i
+\sum\limits_{i=v_1(n)}^{7n+8}
 7^{\mu(n,i)} e_{n,i}u^it
 \\[5pt]
&\quad + \sum\limits_{i=v_2(n)}^{7n+8} 7^{\mu(n,i)}
h_{n,i}u^it^2,
\qquad \qquad\qquad\qquad\qquad\qquad\qquad \text{ if } n\equiv 3 \pmod{7},\\[5pt]
& \sum\limits_{i=v_0(n)}^{7n+8} 7^{\mu(n,i)} d_{n,i}u^i
+\sum\limits_{i=v_1(n)}^{7n+8} 7^{\mu(n,i)}
e_{n,i}u^it+h_{n,v_2(n)}u^{v_2(n)}t^2\\[5pt]
&\quad+ \sum\limits_{i=v_2(n)+1}^{7n+8} 7^{\mu(n,i)}
h_{n,i}u^it^2, \qquad\qquad\qquad\qquad\quad \qquad\qquad\text{ if }
n\equiv 4 \pmod{7}.
\end{cases}
\end{align*}
\end{lemma}

\begin{lemma}
Let $v_0(n)=\left\lceil \frac{n}{7}\right\rceil$,
$v_1(n)=\left\lceil \frac{n+4}{7}\right\rceil$, $v_2(n)=\left\lceil
\frac{n+2}{7}\right\rceil+1$,
$\mu(n,i)=\left\lceil\frac{7i-n}{26}\right\rceil$. Then for
$n\geq 0$, there exist integers $d_{n,i}, e_{n,i}$ and $h_{n,i}$
such that
\begin{align*}
U_7(Au^nt^2)=
\begin{cases}
& 7d_{n,v_0(n)}u^{v_0(n)} +7d_{n,v_0(n)+1}u^{v_0(n)+1}
+\sum\limits_{i=v_0(n)+2}^{7n+8} 7^{\mu(n,i)}
d_{n,i}u^i\\[5pt]
&\quad +e_{n,v_1(n)}u^{v_1(n)}t +e_{n,v_1(n)+1}u^{v_1(n)+1}t
+\sum\limits_{i=v_1(n)+2}^{7n+8} 7^{\mu(n,i)}
e_{n,i}u^it\\[5pt]
&\quad +h_{n,v_2(n)}u^{v_2(n)}t^2+ \sum\limits_{i=v_2(n)+1}^{7n+8}
7^{\mu(n,i)}
h_{n,i}u^it^2,
\qquad\qquad
  \text{ if } n\equiv 1,2,3 \pmod{7},\\[5pt]
& \sum\limits_{i=p_0(n)}^{7n+8} 7^{\mu(n,i)} d_{n,i}u^i
+\sum\limits_{i=v_1(n)}^{7n+8} 7^{\mu(n,i)}
e_{n,i}u^it\\[5pt]
&\quad +h_{n,v_2(n)}u^{v_2(n)}t^2+ \sum\limits_{i=v_2(n)+1}^{7n+8}
7^{\mu(n,i)}
h_{n,i}u^it^2, \qquad
\qquad
\text{ if } n\equiv 4 \pmod{7},\\[5pt]
& d_{n,v_0(n)}u^{v_0(n)}
+d_{n,v_0(n)+1}u^{v_0(n)+1}+d_{n,v_0(n)+2}u^{v_0(n)+2}
\\[5pt]
&\quad+\sum\limits_{i=v_0(n)+3}^{7n+8} 7^{\mu(n,i)} d_{n,i}u^i
+
e_{n,v_1(n)}u^{v_1(n)}t
+ e_{n,v_1(n)+1}u^{v_1(n)+1}t
\\[5pt]
&\quad +\sum\limits_{i=v_1(n)+2}^{7n+8}
 7^{\mu(n,i)}
e_{n,i}u^it+h_{n,v_2(n)}u^{v_2(n)}t^2
\\[5pt]
&\quad +h_{n,v_2(n)+1}u^{v_2(n)+1}t^2+
\sum\limits_{i=v_2(n)+2}^{7n+8} 7^{\mu(n,i)}
h_{n,i}u^it^2, \qquad\
 \text{ if } n\equiv 5 \pmod{7},\\[5pt]
& d_{n,v_0(n)}u^{v_0(n)}
+d_{n,v_0(n)+1}u^{v_0(n)+1}+d_{n,v_0(n)+2}u^{v_0(n)+2}
+\sum\limits_{i=v_0(n)+3}^{7n+8} 7^{\mu(n,i)}
d_{n,i}u^i\\[5pt]
&\quad+ e_{n,v_1(n)}u^{v_1(n)}t+ e_{n,v_1(n)+1}u^{v_1(n)+1}t
+\sum\limits_{i=v_1(n)+2}^{7n+8} 7^{\mu(n,i)}
e_{n,i}u^it\\[5pt]
&\quad +h_{n,v_2(n)}u^{v_2(n)}t^2+ \sum\limits_{i=v_2(n)+1}^{7n+8}
7^{\mu(n,i)} h_{n,i}u^it^2, \qquad\qquad
 \text{ if } n\equiv
0,6 \pmod{7}.
\end{cases}
\end{align*}
\end{lemma}

\subsection{Proof of Theorem \ref{L-expression-7}}

Now, we are ready
 to prove    Theorem \ref{L-expression-7}.

\noindent{\it Proof
 of Theorem \ref{L-expression-7}}.
 We prove the statement by
induction on $\alpha$.
Since $L_1=U_7(A)$, from \eqref{U-7-A} in Appendix B,
we know $L_{2\alpha-1}$ holds for $\alpha=1$.
We assume that $L_{2\alpha-1}$ holds for $\alpha\geq 1$. Applying the
$U_7$-operator to $L_{2\alpha-1}$ yields
\begin{align}
L_{2\alpha}&=U_7(L_{2\alpha-1})\notag\\[5pt]
&=7^{\alpha+1}a_{2\alpha-1,1}U_7(u)+7^{\alpha+1}a_{2\alpha-1,2}U_7(u^2)+
7^{\alpha+1}\sum_{i=3}^{(7^{2\alpha}-1)/6}
7^{\lceil \frac{7i-24}{26}\rceil}a_{2\alpha-1,i}U_7(u^i)\notag\\[5pt]
&\quad
 +7^{\alpha}b_{2\alpha-1,1}U_7(ut)+7^{\alpha}b_{2\alpha-1,2}U_7(u^2t)+
7^{\alpha+1}\sum_{i=3}^{(7^{2\alpha}-1)/6}
7^{\lceil \frac{7i-24}{26}\rceil}b_{2\alpha-1,i}U_7(u^it)\notag\\[5pt]
&\quad + 7^{\alpha}c_{2\alpha-1,2}U_7(u^2t^2)+
7^{\alpha+1}\sum_{i=3}^{(7^{2\alpha}-1)/6} 7^{\lceil
\frac{7i-24}{26}\rceil}c_{2\alpha-1,i} U_7(u^it^2).\label{L-2-comp-7}
\end{align}

Let $S=\{n\colon 3\leq n\leq \frac{7^{2\alpha}-1}{6}\}$, and
$S_i=\{n\in S\colon n\equiv i\pmod{7}\}$ for $0\leq i\leq 6$.
In view of \eqref{U-7-u-t}, \eqref{U-7-u2-t}
and Lemma \ref{U-7-ut}, we obtain that
\begin{align*}
&\quad 7^{\alpha}b_{2\alpha-1,1}U_7(ut)+7^{\alpha}b_{2\alpha-1,2}U_7(u^2t)+
7^{\alpha+1}\sum_{i=3}^{(7^{2\alpha}-1)/6}
7^{\lceil \frac{7i-24}{26}\rceil}b_{2\alpha-1,i}U_7(u^it)\\[5pt]
  & =7^{\alpha}b_{2\alpha-1,1}\bigg(
   7^4 d_{1,1}
   u+7^4d_{1,2} u^2
+\sum_{i=3}^{14}7^{\lceil \frac{7i-11}{4}\rceil}d_{1,i}u^i +7
e_{1,1} ut+ 7^4 e_{1,2}u^2t  +\sum_{i=3}^{14}7^{\lceil
\frac{7i-11}{4}\rceil}e_{1,i}u^it\\[5pt]
&\quad  +7^3 h_{1,2} u^2t^2
 +\sum_{i=3}^{14}7^{\lceil
\frac{7i-11}{4}\rceil}h_{1,i}u^it^2\bigg)
+7^{\alpha}b_{2\alpha-1,2}\bigg(
 7d_{2,1} u+7d_{2,2} u^2+
\sum_{i=3}^{14}7^{\lceil
\frac{7i-12}{4}\rceil}d_{2,i}u^i\\[5pt]
&\quad  +7^2 e_{2,1} ut+7 e_{2,2} u^2t+
 \sum_{i=3}^{14}7^{\lceil
\frac{7i-12}{4}\rceil}e_{2,i}u^it
+7 h_{2,2} u^2t^2+
\sum_{i=3}^{14}7^{\lceil
 \frac{7i-12}{4}\rceil}h_{2,i}u^it^2\bigg)\\[5pt]
&\quad  + 7^{\alpha+1}\sum_{i=3}^{(7^{2\alpha}-1)/6} 7^{\lceil
\frac{7i-24}{26}\rceil}b_{2\alpha-1,i} \left(\sum_{i\in S_0\cup
S_4\cup S_5} \bigg(d_{i,\lceil
 \frac{i}{7}\rceil}u^{\lceil
\frac{i}{7}\rceil} +d_{i,\lceil \frac{i}{7}\rceil+1}u^{\lceil
\frac{i}{7}\rceil+1}\right.\\[5pt]
&\quad  +\sum\limits_{j=\lceil \frac{i}{7}\rceil+2}^{7i} 7^{\lceil
\frac{7j-i+6}{26}\rceil} d_{i,j}u^j +e_{i,\lceil
\frac{i+4}{7}\rceil}u^{\lceil \frac{i+4}{7}\rceil}t
+\sum\limits_{j=\lceil \frac{i+4}{7}\rceil+1}^{7i} 7^{\lceil
\frac{7j-i+6}{26}\rceil} e_{i,j}u^jt \nonumber\\[5pt]
&\quad+h_{i,\lceil \frac{i+2}{7}\rceil+1} u^{\lceil
\frac{i+2}{7}\rceil+1}t^2 + \sum\limits_{j=\lceil
\frac{i+2}{7}\rceil+1}^{7i} 7^{\lceil\frac{7j-i+6}{26}\rceil}
h_{i,j}u^jt^2\bigg)\\[5pt]
&\quad + \sum_{i\in S_1\cup S_2\cup S_3\cup S_6}
 \bigg(
  7d_{i,\lceil
\frac{i}{7}\rceil}u^{\lceil \frac{i}{7}\rceil}
+\sum\limits_{j=\lceil \frac{i}{7}\rceil+1}^{7i}
7^{\lceil\frac{7j-i+6}{26}\rceil} d_{i,j}u^j+ e_{i,\lceil
\frac{i+4}{7}\rceil}u^{\lceil \frac{i+4}{7}\rceil}t
\\[5pt]
&\quad  \left.+\sum\limits_{j=\lceil \frac{i+4}{7}\rceil+1}^{7n}
 7^{\lceil\frac{7j-i+6}{26}\rceil}
e_{i,j}u^jt+ \sum\limits_{i=\lceil \frac{i+2}{7}\rceil+1}^{7i}
7^{\lceil\frac{7j-i+6}{26}\rceil} h_{i,j}u^jt^2\bigg)\right).
\end{align*}
For $i\geq 3$, $\left\lceil\frac{7i-12}{4} \right\rceil\geq
\left\lceil\frac{7i+3}{26} \right\rceil$. For $i\geq 8$,
we have $\lceil
\frac{i}{7}\rceil+1\geq 3$, and
\[\left\lceil \frac{7i-24}{26}\right\rceil
\geq \left\lceil\frac{7(\lceil \frac{i}{7}\rceil+1)+3}{26}
\right\rceil.\] For $i\geq 11$, $\lceil \frac{i+4}{7}\rceil\geq 3$,
and
\[\left\lceil \frac{7i-24}{26}\right\rceil
\geq \left\lceil\frac{7(\lceil \frac{i+4}{7}\rceil)+3}{26}
\right\rceil.\] For $i\geq 6$ and $i\in S_0\cup S_4\cup S_5$,
$\lceil \frac{i+2}{7}\rceil+1\geq 3$, and
\[\left\lceil \frac{7i-24}{26}\right\rceil
\geq \left\lceil\frac{7(\lceil \frac{i+2}{7}\rceil+1)+3}{26}
\right\rceil.\] For $i\geq 15$ and $i\in S_1\cup S_2\cup S_3\cup
S_6$, $\lceil \frac{i}{7}\rceil\geq 3$, and
\[\left\lceil \frac{7i-24}{26}\right\rceil
\geq \left\lceil\frac{7\lceil \frac{i}{7}\rceil+3}{26}
\right\rceil.\]
For $i=3$ and any $j$,
\begin{align*}
\left\lceil\frac{7i-24}{26} \right\rceil +
\left\lceil\frac{7j-i+6}{26} \right\rceil
=\left\lceil\frac{7j+3}{26} \right\rceil.
\end{align*}
For $i\geq 4$ and any $j$,
\begin{align*}
\left\lceil\frac{7i-24}{26} \right\rceil +
\left\lceil\frac{7j-i+6}{26} \right\rceil \geq
\left\lceil\frac{7j+6i-18}{26} \right\rceil
\geq\left\lceil\frac{7j+3}{26} \right\rceil.
\end{align*}
Therefore, there exist integers $a_i,b_i,c_i$ such that
\[7^{\alpha}b_{2\alpha-1,1}U_7(ut)+7^{\alpha}b_{2\alpha-1,2}U_7(u^2t)+
7^{\alpha+1}\sum_{i=3}^{(7^{2\alpha}-1)/6} 7^{\lceil
\frac{7i-24}{26}\rceil}b_{2\alpha-1,i}U_7(u^it)\] can be written as
\begin{align}
&7^{\alpha+1}a_1u+7^{\alpha+1}a_2u^2+
7^{\alpha+1}\sum_{i=3}^{(7^{2\alpha+1}-7)/6} 7^{\lceil
\frac{7i+3}{26}\rceil}a_iu^i
 +7^{\alpha+1}b_1ut+7^{\alpha+1}b_2u^2t\notag\\[5pt]
&\quad + 7^{\alpha+1}\sum_{i=3}^{(7^{2\alpha+1}-7)/6} 7^{\lceil
\frac{7i+3}{26}\rceil}b_iu^it+ 7^{\alpha+1}c_2u^2t^2+
7^{\alpha+1}\sum_{i=3}^{(7^{2\alpha+1}-7)/6} 7^{\lceil
\frac{7i+3}{26}\rceil}c_i u^it^2.\label{L-2-7-form}
\end{align}
In the same vein, each term in \eqref{L-2-comp-7}
can be expressed as the form of \eqref{L-2-7-form}.
Thus, \eqref{L-2n} holds for $\alpha$. Therefore,
\begin{align}
L_{2\alpha+1}&=U_7(AL_{2\alpha})\notag\\[5pt]
&=7^{\alpha+1}a_{2\alpha,1}U_7(Au)
+7^{\alpha+1}a_{2\alpha,2}U_7(Au^2) +
7^{\alpha+1}\sum_{i=3}^{(7^{2\alpha+1}-7)/6}
7^{\lceil \frac{7i+3}{26}\rceil}a_{2\alpha,i}U_7(Au^i)\notag\\[5pt]
&\quad
 +7^{\alpha+1}b_{2\alpha,1}U_7(Aut)+7^{\alpha+1}b_{2\alpha,2}U_7(Au^2t) +
7^{\alpha+1}\sum_{i=3}^{(7^{2\alpha+1}-7)/6}
7^{\lceil \frac{7i+3}{26}\rceil}b_{2\alpha,i}U_7(Au^it)\notag\\[5pt]
&\quad+ 7^{\alpha+1}c_{2\alpha,2}U_7(Au^2t^2)+
7^{\alpha+1}\sum_{i=3}^{(7^{2\alpha+1}-7)/6} 7^{\lceil
\frac{7i+3}{26}\rceil}c_{2\alpha,i} U_7(Au^it^2).\label{L-odd-prove}
\end{align}
Then, utilizing the expressions of $U_7(Au)$ and $U_7(Au^2)$
in our Mathematica supplement,
and Lemma \ref{U-7-A-u}, we have
\begin{align*}
&\quad7^{\alpha+1}a_{2\alpha,1}U_7(Au) +7^{\alpha+1}a_{2\alpha,2}U_7(Au^2) +
7^{\alpha+1}\sum_{i=3}^{(7^{2\alpha+1}-7)/6}
7^{\lceil \frac{7i+3}{26}\rceil}a_{2\alpha,i}U_7(Au^i)\\[5pt]
&=7^{\alpha+1}a_{2\alpha,1}(7 d_{1,1} u+7d_{1,2} u^2
+\sum_{i=3}^{15}7^{\lceil\frac{7i-17}{4}\rceil}d_{1,i}u^i +e_{1,1}
ut+e_{1,2} u^2t+\sum_{i=3}^{15}
7^{\lceil\frac{7i-17}{4}\rceil}e_{1,i}u^it\\[5pt]
&\quad +h_{1,2}u^2t^2+\sum_{i=3}^{15}7^{\lceil\frac{7i-17}{4}\rceil}
h_{1,i}u^it^2) +7^{\alpha+1}a_{2\alpha,2}(7^2 d_{2,1} u+ 7^2d_{2,2}
u^2
+\sum_{i=3}^{15}7^{\lceil\frac{7i-18}{4}\rceil}d_{2,i}u^i\\[5pt]
&\quad +7^2e_{2,2} u^2t+\sum_{i=3}^{15}
7^{\lceil\frac{7i-18}{4}\rceil}e_{2,i}u^it
+7h_{2,2}u^2t^2+\sum_{i=3}^{15}7^{\lceil\frac{7i-18}{4}\rceil}
h_{2,i}u^it^2)\\[5pt]
&\quad+ 7^{\alpha+1}\sum_{i=3}^{(7^{2\alpha+1}-7)/6} 7^{\lceil
\frac{7i+3}{26}\rceil}a_{2\alpha,i} (d_{i,\lceil
\frac{i+1}{7}\rceil}u^{\lceil \frac{i+1}{7}\rceil} +d_{i,\lceil
\frac{i+1}{7}\rceil+1}u^{\lceil \frac{i+1}{7}\rceil+1}+ d_{i,\lceil
\frac{i+1}{7}\rceil+2}
u^{\lceil \frac{i+1}{7}\rceil+2}\\[5pt]
&\quad +\sum\limits_{j=\lceil \frac{i+1}{7}\rceil+3}^{7i+8}
 7^{\lceil\frac{7j-i}{26}\rceil}
d_{i,j}u^j+ e_{i,\lceil \frac{i+6}{7}\rceil}u^{\lceil
\frac{i+6}{7}\rceil}t +e_{i,\lceil \frac{i+6}{7}\rceil+1}u^{\lceil
\frac{i+6}{7}\rceil+1}t +\sum\limits_{j=\lceil
\frac{i+6}{7}\rceil+2}^{7i+8} 7^{\lceil\frac{7j-i}{26}\rceil}
e_{i,j}u^jt\\[5pt]
&\quad +h_{i,\lceil \frac{i+3}{7}\rceil+1} u^{\lceil
\frac{i+3}{7}\rceil+1}t + \sum\limits_{j=\lceil
\frac{i+3}{7}\rceil+2}^{7i+8}
 7^{\lceil\frac{7j-i}{26}\rceil}
h_{i,j}u^jt^2).
\end{align*}
For $i\geq 3$,
\[\left\lceil \frac{7i-22}{4}\right\rceil\geq
\left\lceil \frac{7i-24}{26}\right\rceil.\] For $i\geq 3$,
$\lceil\frac{i+1}{7} \rceil\geq 1$, $\lceil\frac{i+6}{7} \rceil\geq
2$, and
\[\left\lceil \frac{7i-23}{26}\right\rceil\geq
\left\lceil \frac{7(\lceil\frac{i+1}{7}
\rceil+2)-24}{26}\right\rceil\geq \left\lceil
\frac{7(\lceil\frac{i+6}{7} \rceil+1)-24}{26}\right\rceil.\]
So, for
$i\geq 5$,
\[\left\lceil \frac{7i-23}{26}\right\rceil\geq
\left\lceil \frac{7(\lceil\frac{i+3}{7}
\rceil+1)-24}{26}\right\rceil.\]
For $i\geq 3$ and any $j$,
\begin{align*}
\left\lceil \frac{7i-23}{26}\right\rceil+ \left\lceil
\frac{7j-i}{26}\right\rceil \geq \left\lceil
\frac{7j+6i-23}{26}\right\rceil \geq \left\lceil
\frac{7j-24}{26}\right\rceil.
\end{align*}
Therefore, there exist integers $a_i,b_i,c_i$ such that
\[7^{\alpha+1}a_{2\alpha,1}U_7(Au)
+7^{\alpha+1}a_{2\alpha,2}U_7(Au^2) +
7^{\alpha+1}\sum_{i=3}^{(7^{2\alpha+1}-7)/6} 7^{\lceil
\frac{7i+3}{26}\rceil}a_{2\alpha,i}U_7(Au^i)\] can be expressed as
\begin{align*}
&7^{\alpha+2}a_1u+7^{\alpha+2}a_2u^2+
7^{\alpha+2}\sum_{i=3}^{(7^{2\alpha+2}-1)/6} 7^{\lceil
\frac{7i-24}{26}\rceil}a_iu^i
 +7^{\alpha+1}b_1ut+7^{\alpha+1}b_2u^2t\\[5pt]
&\quad+ 7^{\alpha+2}\sum_{i=3}^{(7^{2\alpha+2}-1)/6} 7^{\lceil
\frac{7i-24}{26}\rceil}b_iu^it + 7^{\alpha+1}c_2u^2t^2+
7^{\alpha+2}\sum_{i=3}^{(7^{2\alpha+2}-1)/6} 7^{\lceil
\frac{7i-24}{26}\rceil}c_i u^it^2.
\end{align*}
In the similar way, we can obtain the other terms of
\eqref{L-odd-prove} can be expressed as the above form. Hence,
\eqref{L-2n-1} holds for $\alpha+1$.
By induction, the proof is complete.
\qed

\section{Concluding Remarks}

Over the last century,
 a large variety of infinite congruence
families have been discovered and studied, exhibiting a great
variety with respect to their difficulty; see
 for example \cite{Paule-Radu,Smoot-1}.
 In this paper, we investigate
 congruence properties
  for $PDN1(N)$ which
  denotes  the number of partition
   diamonds with $(n+1)$ copies of $n$
    where summing the parts
     at the links gives $N$. The partition
      function  $PDN1(N)$ has received
       a lot of attention
 in recent years. Motivated by Andrews and Paule's
  work \cite{Andrews-2024}, we prove
    three congruences modulo
 arbitrary powers of
 $5$  for $PDN1(n)$
  and  two congruences modulo
 arbitrary powers of
 $7$  for $PDN1(n)$.
  In particular, we  give an answer
  of Andrews and Paule's open problem
  on  determining
               infinite families of
                congruences similar to
                Ramanujan's classical $
                 p(5^kn +d_k) \equiv
                  0 \pmod {5^k}$, where
             $24d_k\equiv 1 \pmod {5^k}$
and $k\geq 1$.

\section{Acknowledgements}

Julia Q.D. Du was supported by the National Natural Science
Foundation of China (No. 12201177), the Hebei Natural Science
Foundation (Nos.~A2024205012, A2023205045), the Science and
Technology Project of Hebei Education Department (No. BJK2023092),
the China Scholarship Council, the Program for
Foreign Experts of Hebei Province, and
the Special Project on Science and Technology Research and Development
Platforms, Hebei Province (22567610H). Olivia X.M. Yao was
 supported by the Natural Science
Foundation of Jiangsu Province of China (No. BK20221383).

\noindent{\bf Author Contributions.}
The authors contributed equally to the preparation of this article. All authors read and
 approved the final manuscript.

\noindent{\bf Competing Interests.}
The authors declare that they have
no conflict of interest.

\noindent{\bf Data Availability Statements.} 
Data sharing not applicable to this
article as no datasets were generated or analysed during the current
study.

\section{Appendix}\label{appendix}

\subsection*{Appendix A}

\begin{align}
U_5(A)
&=-(3 u + u t + u^2 t)+
8675 u + 3277600 u^2 + 295664000 u^3 + 10852800000 u^4\notag\\[5pt]
 &\quad +
 195776000000 u^5 + 1840000000000 u^6 + 8640000000000 u^7 +
 16000000000000 u^8\notag\\[5pt]
 &\quad+
 150 ut + 193825 u^2t + 28016000 u^3t + 1343680000 u^4t +
 28940800000 u^5t \notag\\[5pt]
 &\quad + 309120000000 u^6t+ 1600000000000 u^7t +
 3200000000000 u^8t,\label{U-A}\\[5pt]
U_5(Au)&=-2 (3 u + u t + u^2 t)
+7650 u + 30 t u + 19392110 u^2 + 517010 t u^2 + 7890080750 u^3\notag\\[5pt]
&\quad +
 410659050 t u^3 + 1168362087500 u^4 + 88030965625 t u^4+
 85209875359375 u^5 \notag\\[5pt]
&\quad + 8186499490625 t u^5 + 3551479600000000 u^6 +
 406471523203125 t u^6 \notag\\[5pt]
&\quad+ 91749032000000000 u^7 +
 12000545600000000 t u^7 + 1534502400000000000 u^8\notag\\[5pt]
&\quad +
 223135360000000000 t u^8 + 16902000000000000000 u^9 +
 2679792000000000000 t u^9 \notag\\[5pt]
&\quad+ 121713600000000000000 u^{10} +
 20742080000000000000 t u^{10} + 551200000000000000000 u^{11} \notag\\[5pt]
&\quad+
 99872000000000000000 t u^{11} + 1424000000000000000000 u^{12} +
 272000000000000000000 t u^{12} \notag\\[5pt]
&\quad+ 1600000000000000000000 u^{13} +
 320000000000000000000 t u^{13}\label{A-u-f}\\[5pt]
 &=2 (7 u - u t + 4 u^2 t - 15 u^2)+
 7630 u + 30 t u + 19392140 u^2 + 517000 t u^2+ 7890080750 u^3 \notag\\[5pt]
&\quad +
 410659050 t u^3 + 1168362087500 u^4 + 88030965625 t u^4 +
 85209875359375 u^5 + 8186499490625 t u^5  \notag\\[5pt]
&\quad+ 3551479600000000 u^6 +
 406471523203125 t u^6 + 91749032000000000 u^7 +
 12000545600000000 t u^7  \notag\\[5pt]
&\quad+ 1534502400000000000 u^8 +
 223135360000000000 t u^8 + 16902000000000000000 u^9  \notag\\[5pt]
&\quad+
 2679792000000000000 t u^9 + 121713600000000000000 u^{10} +
 20742080000000000000 t u^{10} \notag\\[5pt]
&\quad + 551200000000000000000 u^{11} +
 99872000000000000000 t u^{11} + 1424000000000000000000 u^{12} \notag\\[5pt]
&\quad +
 272000000000000000000 t u^{12} + 1600000000000000000000 u^{13}\notag\\[5pt]
&\quad +
 320000000000000000000 t u^{13},\label{A-u-s}\\[5pt]
U_5(Au^2)&=(3 u + u t + u^2 t)+2210 u + 32838410 u^2 + 419980 t u^2 + 45424873625 u^3 +
 1394567975 t u^3 \notag\\[5pt]
&\quad+ 19218691600000 u^4 + 939986728125 t u^4 +
 3711979441765625 u^5 + 245455242090625 t u^5\notag\\[5pt]
&\quad +
 397914022469531250 u^6 + 32647890206421875 t u^6 +
 26478069076474609375 u^7 \notag\\[5pt]
&\quad+ 2559835121300781250 t u^7 +
 1172753998329394531250 u^8 + 129128836223105468750 t u^8 \notag\\[5pt]
&\quad+
 36174571382399658203125 u^9 + 4429647038721093750000 t u^9 +
 800131993455999755859375 u^{10}\notag\\[5pt]
&\quad + 107077436156801757812500 t u^{10} +
 12914193040000000000000000 u^{11}\notag\\[5pt]
&\quad + 1863895822720001220703125 t u^{11} +
 153270802400000000000000000 u^{12}\notag\\[5pt]
&\quad +
 23612522080000000000000000 t u^{12} +
 1335232640000000000000000000 u^{13}\notag\\[5pt]
&\quad +
 217763392000000000000000000 t u^{13} +
 8429832000000000000000000000 u^{14}\notag\\[5pt]
&\quad +
 1445697600000000000000000000 t u^{14} +
 37504960000000000000000000000 u^{15}\notag\\[5pt]
&\quad +
 6726208000000000000000000000 t u^{15} +
 111440000000000000000000000000 u^{16}\notag\\[5pt]
&\quad +
 20803200000000000000000000000 t u^{16} +
 198400000000000000000000000000 u^{17}\notag\\[5pt]
&\quad +
 38400000000000000000000000000 t u^{17} +
 160000000000000000000000000000 u^{18}\notag\\[5pt]
&\quad +
 32000000000000000000000000000 t u^{18}\label{A-u-2-f}\\[5pt]
 &=-(7 u - u t + 4 u^2 t - 15 u^2)+
2220 u + 32838395 u^2 + 419985 t u^2 + 45424873625 u^3 +
 1394567975 t u^3 \notag\\[5pt]
&\quad+ 19218691600000 u^4 + 939986728125 t u^4 +
 3711979441765625 u^5 + 245455242090625 t u^5 \notag\\[5pt]
&\quad+
 397914022469531250 u^6 + 32647890206421875 t u^6 +
 26478069076474609375 u^7 \notag\\[5pt]
&\quad+ 2559835121300781250 t u^7 +
 1172753998329394531250 u^8 + 129128836223105468750 t u^8 \notag\\[5pt]
&\quad+
 36174571382399658203125 u^9 + 4429647038721093750000 t u^9 +
 800131993455999755859375 u^{10}\notag\\[5pt]
&\quad + 107077436156801757812500 t u^{10} +
 12914193040000000000000000 u^{11} \notag\\[5pt]
&\quad+ 1863895822720001220703125 t u^{11}
  +
 153270802400000000000000000 u^{12} \notag\\[5pt]
&\quad+
 23612522080000000000000000 t u^{12} +
 1335232640000000000000000000 u^{13} \notag\\[5pt]
&\quad+
 217763392000000000000000000 t u^{13} +
 8429832000000000000000000000 u^{14} \notag\\[5pt]
&\quad+
 1445697600000000000000000000 t u^{14} +
 37504960000000000000000000000 u^{15} \notag\\[5pt]
&\quad+
 6726208000000000000000000000 t u^{15} +
 111440000000000000000000000000 u^{16}\notag\\[5pt]
&\quad +
 20803200000000000000000000000 t u^{16} +
 198400000000000000000000000000 u^{17} \notag\\[5pt]
&\quad+
 38400000000000000000000000000 t u^{17} +
 160000000000000000000000000000 u^{18} \notag\\[5pt]
&\quad+
 32000000000000000000000000000 t u^{18},\label{A-u-2-s}\\[5pt]
U_5(Aut)&=-2 (3 u + u t + u^2 t)+1735 u + 30 t u + 655370 u^2 + 38690 t u^2 + 59132250 u^3 +
 5602350 t u^3 \notag\\[5pt]
&\quad+ 2170562500 u^4 + 268731875 t u^4 + 39155203125 u^5 +
 5788146875 t u^5 + 368000000000 u^6\notag\\[5pt]
&\quad + 61823984375 t u^6 +
 1728000000000 u^7 + 320000000000 t u^7 + 3200000000000 u^8 \notag\\[5pt]
&\quad+
 640000000000 t u^8\label{A-u-t-f}\\[5pt]
&=2 (7 u - u t + 4 u^2 t - 15 u^2)+
1715 u + 30 t u + 655400 u^2 + 38680 t u^2 + 59132250 u^3 +
 5602350 t u^3 \notag\\[5pt]
&\quad + 2170562500 u^4 + 268731875 t u^4 + 39155203125 u^5 +
 5788146875 t u^5 + 368000000000 u^6 \notag\\[5pt]
&\quad + 61823984375 t u^6 +
 1728000000000 u^7 + 320000000000 t u^7 + 3200000000000 u^8\notag\\[5pt]
&\quad  +
 640000000000 t u^8,\label{A-u-t-s}\\[5pt]
U_5(Au^2t)&=-(3 u + u t + u^2 t)+
1475 u + 5 t u + 3853385 u^2 + 101685 t u^2 + 1575693125 u^3 +
 81896675 t u^3 \notag\\[5pt]
&\quad+ 233586412500 u^4 + 17595116250 t u^4 +
 17040418500000 u^5 + 1637062768750 t u^5 \notag\\[5pt]
&\quad+ 710281272343750 u^6 +
 81291771250000 t u^6 + 18349737626953125 u^7 +
 2400095896093750 t u^7 \notag\\[5pt]
&\quad+ 306900353027343750 u^8 +
 44627044472656250 t u^8 + 3380400001708984375 u^9\notag\\[5pt]
&\quad +
 535958394531250000 t u^9 + 24342720001220703125 u^{10} +
 4148415991210937500 t u^{10}\notag\\[5pt]
&\quad + 110240000000000000000 u^{11} +
 19974399993896484375 t u^{11} + 284800000000000000000 u^{12}\notag\\[5pt]
&\quad +
 54400000000000000000 t u^{12} + 320000000000000000000 u^{13} \notag\\[5pt]
&\quad+
 64000000000000000000 t u^{13}\label{A-u-2-t-f}\\[5pt]
 &=(7 u - u t + 4 u^2 t - 15 u^2)+
 1465 u + 5 t u + 3853400 u^2 + 101680 t u^2 + 1575693125 u^3 +
 81896675 t u^3 \notag\\[5pt]
&\quad+ 233586412500 u^4 + 17595116250 t u^4 +
 17040418500000 u^5 + 1637062768750 t u^5 \notag\\[5pt]
&\quad+ 710281272343750 u^6 +
 81291771250000 t u^6 + 18349737626953125 u^7 +
 2400095896093750 t u^7 \notag\\[5pt]
&\quad+ 306900353027343750 u^8 +
 44627044472656250 t u^8 + 3380400001708984375 u^9\notag\\[5pt]
&\quad +
 535958394531250000 t u^9 + 24342720001220703125 u^{10} +
 4148415991210937500 t u^{10}\notag\\[5pt]
&\quad + 110240000000000000000 u^{11} +
 19974399993896484375 t u^{11} + 284800000000000000000 u^{12}\notag\\[5pt]
&\quad +
 54400000000000000000 t u^{12} + 320000000000000000000 u^{13}\notag\\[5pt]
&\quad +
 64000000000000000000 t u^{13},\label{A-u-2-t-s}\\[5pt]
U_5(3 u &+ u t + u^2 t)= 10 (u t - u^2 t)
+525 u + 25 t u + 58075 u^2 + 9300 t u^2 + 1485625 u^3 +
 420000 t u^3 \notag\\[5pt]
&\quad+ 10125000 u^4 + 7443750 t u^4 - 44062500 u^5 +
 78640625 t u^5 - 734375000 u^6 + 667187500 t u^6 \notag\\[5pt]
&\quad- 2978515625 u^7 +
 4218750000 t u^7 - 6591796875 u^8 + 16162109375 t u^8 -
 6103515625 u^9 \notag\\[5pt]
&\quad+ 34179687500 t u^9 + 30517578125 t u^{10},\label{U-3-u}\\[5pt]
U_5(7 u &- u t + 4 u^2 t - 15 u^2)=-10 (u - u^2 t)
-5225 u + 25 t u - 6165075 u^2 - 217800 t u^2 - 1205398125 u^3\notag\\[5pt]
&\quad -
 84020000 t u^3 - 89462671875 u^4 - 8730771875 t u^4 -
 3280968437500 u^5 - 397291531250 t u^5\notag\\[5pt]
&\quad - 66798640625000 u^6 -
 9437229687500 t u^6 - 790820849609375 u^7 - 125610468750000 t u^7\notag\\[5pt]
&\quad -
 5409046142578125 u^8 - 942486865234375 t u^8 -
 19800042724609375 u^9 - 3719760742187500 t u^9\notag\\[5pt]
&\quad -
 30000000000000000 u^{10} - 5999786376953125 t u^{10},\label{U-7-u}\\[5pt]
U_5(Au&-Au^2t)=- (3 u + u t + u^2 t)+
6175 u + 25 t u + 15538725 u^2 + 415325 t u^2 + 6314387625 u^3\notag\\[5pt]
&\quad +
 328762375 t u^3 + 934775675000 u^4 + 70435849375 t u^4 +
 68169456859375 u^5 \notag\\[5pt]
&\quad+ 6549436721875 t u^5 + 2841198327656250 u^6 +
 325179751953125 t u^6 + 73399294373046875 u^7\notag\\[5pt]
&\quad +
 9600449703906250 t u^7 + 1227602046972656250 u^8 +
 178508315527343750 t u^8 \notag\\[5pt]
&\quad+ 13521599998291015625 u^9 +
 2143833605468750000 t u^9 + 97370879998779296875 u^{10}\notag\\[5pt]
&\quad +
 16593664008789062500 t u^{10} + 440960000000000000000 u^{11} +
 79897600006103515625 t u^{11}\notag\\[5pt]
&\quad + 1139200000000000000000 u^{12} +
 217600000000000000000 t u^{12}\notag\\[5pt]
&\quad + 1280000000000000000000 u^{13} +
 256000000000000000000 t u^{13},\label{U-u-u2}\\[5pt]
U_5(Aut& - Au^2t)=(7 u - u t + 4 u^2 t - 15 u^2)+
250 u + 25 t u - 3198000 u^2 - 63000 t u^2 - 1516560875 u^3\notag\\[5pt]
&\quad -
 76294325 t u^3 - 231415850000 u^4 - 17326384375 t u^4 -
 17001263296875 u^5 - 1631274621875 t u^5\notag\\[5pt]
&\quad - 709913272343750 u^6 -
 81229947265625 t u^6 - 18348009626953125 u^7 -
 2399775896093750 t u^7\notag\\[5pt]
&\quad - 306897153027343750 u^8 -
 44626404472656250 t u^8 - 3380400001708984375 u^9\notag\\[5pt]
&\quad -
 535958394531250000 t u^9 - 24342720001220703125 u^{10} -
 4148415991210937500 t u^{10}\notag\\[5pt]
&\quad - 110240000000000000000 u^{11} -
 19974399993896484375 t u^{11} - 284800000000000000000 u^{12}\notag\\[5pt]
&\quad -
 54400000000000000000 t u^{12} - 320000000000000000000 u^{13}\notag\\[5pt]
&\quad -
 64000000000000000000 t u^{13}.\label{U-ut-u2t}
\end{align}

\subsection*{Appendix B}

\begin{align*}
\sigma _1(q)&=(1766604 u+3524494638 u^2+477447678477 u^3+3979690016743 u^4-123260426116415 u^5\\[5pt]
&\quad+157724494175920 u^6+2517690595327040 u^7)+(434 u + 43880284 u^2 + 45959089645 u^3 \\[5pt]
&\quad+ 5000367092931 u^4 +
 61462358112176 u^5 - 647031711964743 u^6\\[5pt]
&\quad - 2897179143843600 u^7)t
+(51303 u^2 + 513960461 u^3 + 217979361306 u^4 + 14148972748316 u^5\\[5pt]
&\quad +
 183847805203440 u^6 + 379502389803761 u^7)t^2,\\[5pt]
\sigma_2(q)&=(-811706 u - 2296108591 u^2 - 341824185024 u^3 -
3114592566400 u^4 +
 90823431738368 u^5 \\[5pt]
&\quad- 107933308399616 u^6 - 1878998238724096 u^7)+
 (21 u - 23203852 u^2 - 30503994304 u^3\\[5pt]
 &\quad - 3579966833184 u^4 -
 46238421274304 u^5 + 471423104470528 u^6 + 2162467800600576 u^7)t\\[5pt]
 &\quad+
 (-15827 u^2 - 306624507 u^3 - 149208932495 u^4 - 10184367568704 u^5 -
 135557127949824 u^6\\[5pt]
 &\quad - 283469561876480 u^7)t^2,\\[5pt]
\sigma_3(q)&= (161637 u + 630993580 u^2 + 101953152959 u^3 +
1007777300320 u^4 -
 27865427301376 u^5 \\[5pt]
 &\quad+ 30599010287104 u^6 + 584294403051520 u^7)+
 (5347664 u^2 + 8521047815 u^3\\[5pt]
 &\quad + 1068550874496 u^4 +
 14462166151216 u^5 - 143116488441920 u^6 - 672517072819200 u^7)t\\[5pt]
 &\quad+
 (1617 u^2 + 78480458 u^3 + 42812672784 u^4 + 3057610920336 u^5 +
 41650087685696 u^6\\[5pt]
 &\quad + 88222669767680 u^7)t^2,\\[5pt]
\sigma_4(q)&=(-18648 u - 93561391 u^2 - 16221979886 u^3 -
172775304527 u^4 +
 4556953306136 u^5 \\[5pt]
 &\quad- 4595976067648 u^6 - 96900309416960 u^7)+
 (-695296 u^2 - 1281364994 u^3 - 170232960058 u^4\\[5pt]
 &\quad - 2407848670024 u^5 +
 23172701402088 u^6 + 111543826325120 u^7)t+
 (-56 u^2 - 11028136 u^3\\[5pt]
 &\quad - 6588091895 u^4 - 490087145128 u^5 -
 6825754976040 u^6 - 14643516908160 u^7)t^2,\\[5pt]
\sigma_5(q)&=(1344 u + 7886060 u^2 + 1452706871 u^3 + 16569113722
u^4 -
 418966089297 u^5 + 385273180432 u^6\\[5pt]
 &\quad + 9039207968000 u^7)+
 (53704 u^2 + 109301066 u^3 + 15269092272 u^4 + 225121479149 u^5\\[5pt]
 &\quad -
 2110568588513 u^6 - 10406388173160 u^7)t+(894201 u^3 + 573152314 u^4\\[5pt]
 &\quad + 44230729795 u^5 + 629291441961 u^6 +
 1367180205160 u^7)t^2,\\[5pt]
\sigma_6(q)&=(-56 u - 357798 u^2 - 69439664 u^3 - 843418478 u^4 +
20534456460 u^5 -
 17063810960 u^6\\[5pt]
 &\quad - 449700596408 u^7)+
 (-2303 u^2 - 5009515 u^3 - 731186134 u^4 - 11208890826 u^5\\[5pt]
 &\quad  +
 102526985785 u^6+ 517777131417 u^7)t+
 (-39445 u^3 - 26715927 u^4 \\[5pt]
 &\quad- 2131093986 u^5 - 30945451768 u^6 -
 68076535009 u^7)t^2,\\[5pt]
\sigma_7(q)&=(u + 6818 u^2 + 1384397 u^3 + 17815420 u^4 - 419183387
u^5 +
 311299254 u^6 + 9321683217 u^7)\\[5pt]
 &\quad+
 (42 u^2 + 96285 u^3 + 14605283 u^4 + 232239126 u^5 - 2075328360 u^6 -
 10734059462 u^7)t\\[5pt]
 &\quad+
 (735 u^3 + 521017 u^4 + 42824236 u^5 + 634128110 u^6 + 1412376245 u^7)t^2.
\end{align*}

\begin{align}
U_7(A)&=
2647568 u + 10203888576 u^2 + 2817650339840 u^3 +
 108003640500224 u^4 - 560421350866944 u^5\notag\\[5pt]
&\quad - 14565886187274240 u^6 +
 59747119224848384 u^7 + 203783375486451712 u^8\notag\\[5pt]
&\quad+
 (371 u + 78499344 u^2 + 148515453184 u^3 + 31424911085568 u^4 +
 1179224675844096 u^5 \notag\\[5pt]
&\quad+ 1982920774123520 u^6 -
 107371327927091200 u^7 - 233403052155994112 u^8)t
 \notag\\[5pt]
&\quad
+(59185 u^2 + 1147512576 u^3 + 871469813760 u^4 + 107224172396544 u^5\notag\\[5pt]
&\quad +
 3137915454750720 u^6 + 22743680299827200 u^7 + 29619676669542400 u^8)t^2,\label{U-7-A}\\[5pt]
U_7(ut)&=
2028845 u + 574 t u + 4028076066 u^2 + 50260133 t u^2 +
 60368 t^2 u^2 + 545647905803 u^3\notag\\[5pt]
&\quad + 52525032315 t u^3 +
 587675963 t^2 u^3 + 4548338693409 u^4 + 5714686879013 t u^4\notag\\[5pt]
&\quad +
 249121051606 t^2 u^4 - 140870389381897 u^5 + 70242854987984 t u^5 +
 16170243745796 t^2 u^5\notag\\[5pt]
&\quad + 180264404901840 u^6 -
 739466249109441 t u^6 + 210111869612176 t^2 u^6 +
 2877342601957824 u^7 \notag\\[5pt]
&\quad- 3311057359074416 t u^7 +
 433716734443335 t^2 u^7,\label{U-7-u-t}\\[5pt]
 U_7(u^2 t)&=54794187 u + 147 t u + 9410867661755 u^2 + 5467251020 t u^2 +
 235165 t^2 u^2 + 72534495322245383 u^3 \notag\\[5pt]
&\quad+ 295146766935111 t u^3 +
 286211620373 t^2 u^3 + 97274018997437919364 u^4 \notag\\[5pt]
&\quad+
 1392226546484060392 t u^4 + 5404244180728169 t^2 u^4 +
 37138586840405605841853 u^5\notag\\[5pt]
&\quad + 1425047059370490952571 t u^5 +
 13639843118966429869 t^2 u^5 + 5018460780942093521877505 u^6\notag\\[5pt]
&\quad +
 462638465634577854814059 t u^6 + 8980480326202590208822 t^2 u^6 \notag\\[5pt]
&\quad+
 255616325942673945513397238 u^7 + 57205701288543441686138044 t u^7\notag\\[5pt]
&\quad +
 2094619590855642617865819 t^2 u^7 +
 4128666559770447770399082661 u^8 \notag\\[5pt]
&\quad+
 2891568889508876423502068350 t u^8 +
 203110006153574678274539409 t^2 u^8 \notag\\[5pt]
&\quad-
 15553360645704583425971672768 u^9 +
 56464577399786426567233403076 t u^9 \notag\\[5pt]
&\quad+
 8897230756527228307709916700 t^2 u^9 -
 626347135495368324670585070209 u^{10} \notag\\[5pt]
&\quad+
 217190999304984710675108330472 t u^{10} +
 181720475986679847777219957612 t^2 u^{10}\notag\\[5pt]
&\quad -
 764344804181577940902407076642 u^{11} -
 4045056394786993046189713820122 t u^{11}\notag\\[5pt]
&\quad +
 1706834305729161181616735590024 t^2 u^{11} +
 20081584017396471682549320033031 u^{12}\notag\\[5pt]
&\quad -
 34045489011853008713079465422139 t u^{12} +
 6877251570376423176209512683517 t^2 u^{12}\notag\\[5pt]
&\quad +
 56018783609819956029469538070800 u^{13} -
 69219297342053944591140361729457 t u^{13}\notag\\[5pt]
&\quad +
 10053109412154871774238735547885 t^2 u^{13} +
 25486298396574955116884857981376 u^{14}\notag\\[5pt]
&\quad -
 28940766552705567711105422852720 t u^{14} +
 3454468156157981341560645787687 t^2 u^{14}.\label{U-7-u2-t}
\end{align}

\end{document}